\documentclass[reqno,a4paper]{amsart}
\usepackage{graphicx} 
\usepackage[utf8]{inputenc}
\usepackage[T1]{fontenc}
\usepackage{amsmath, amssymb, amscd, amsthm,epsfig,epic,eepic}
\usepackage[english]{babel}
\usepackage{mathrsfs}
\usepackage{tikz-cd}
\usepackage{quiver}
\usepackage{verbatim}
\usepackage{color}
\usepackage[normalem]{ulem}
\usepackage[all]{xy}
\usepackage[justify]{ragged2e}
\usepackage{csquotes}
\usepackage[backend=bibtex,style=alphabetic]{biblatex}
\usepackage{hyperref, cleveref}
\hypersetup{
    colorlinks=true,
    linkcolor=teal,
    citecolor=teal,
    filecolor=teal,
    urlcolor=teal,
    linktocpage=true
}
\usepackage[left=3cm,right=3cm,top=3cm,bottom=3cm]{geometry}

\theoremstyle{plain}

\newtheorem{theorem}{Theorem}[section]
\newtheorem{lemma}[theorem]{Lemma}
\newtheorem{proposition}[theorem]{Proposition}
\newtheorem{corollary}[theorem]{Corollary}
\newtheorem{theorem*}{Theorem}
\newtheorem{proposition*}[theorem*]{Proposition}

\newtheorem{theoremA}{Theorem}

\newtheorem{propositionA}[theoremA]{Proposition}

\newtheorem{corollaryA}[theoremA]{Corollary}

\theoremstyle{definition}
\newtheorem{definition}[theorem]{Definition}
\newtheorem{remark}[theorem]{Remark}

\newtheorem{example}[theorem]{Example}
\numberwithin{equation}{section}

\newcommand{\blfootnote}[1]{%
  \begingroup
  \renewcommand{\thefootnote}{}%
  \begin{NoHyper}
    \footnote{#1}%
  \end{NoHyper}
  \addtocounter{footnote}{-1}%
  \endgroup
}
 
\makeatletter
     \let\oldfootnote\footnote
     \def\footnote{\@ifstar\footnote@star\footnote@nostar}    \def\footnote@star#1{{\let\thefootnote\relax\footnotetext{#1}}}
     \def\footnote@nostar{\oldfootnote}
     \makeatother

\def\Spec{\operatorname{Spec}}

\renewcommand{\P}{\mathbb{P}}
\renewcommand{\O}{\mathcal{O}}

\DeclareMathOperator{\muu}{\boldsymbol{\mu}}
\DeclareMathOperator{\aalpha}{\boldsymbol{\alpha}}
\DeclareMathOperator{\Pic}{\underline{Pic}}
\DeclareMathOperator{\Aut}{\underline{Aut}}
\DeclareMathOperator{\NS}{\underline{NS}}

\title{Generalizations of Igusa surfaces}

\author{Pascal Fong}
   \address{Leibniz Universität Hannover, Institut für Algebraische Geometrie, Welfengarten 1, 30167 Hannover}
\email{fong@math.uni-hannover.de} 

\author{Matilde Maccan}
   \address{Institut Camille Jordan, UMR CNRS 5208, Université Lyon 1, Villeurbanne, France}
\email{maccan@math.univ-lyon1.fr} 

\author{Ludvig Modin}
\address{Leibniz Universität Hannover, Institut für Algebraische Geometrie, Welfengarten 1, 30167 Hannover}
\email{modin@math.uni-hannover.de}

\date{\today}

\begin{document}

\begin{abstract}
   We give a criterion for the contracted product \(E \times^G X\) to have a singular Picard scheme, where \(E\) is an elliptic curve, \(G \subset E\) is a finite subgroup scheme, and \(X\) is a projective \(G\)-normal curve, in the sense of Brion. We show that every smooth projective surface with a singular Picard scheme and a positive-dimensional automorphism group scheme arises in this way. We also illustrate the criterion by constructing two new families of surfaces with singular Picard schemes, arising from \(G\)-normal curves \(X\) for \(G = \mathbb{Z}/p\mathbb{Z}\) and \(G = \aalpha_p\), respectively. As a special case of \(G = \mathbb{Z}/2\mathbb{Z}\), we recover the Igusa’s classical example. Finally, we compute all infinitesimal group schemes which appear as \(\Pic^0_S/(\Pic^0_S)_{\rm red}\) for all (quasi)-bielliptic surfaces $S$. This completes the description of the singular Picard schemes of smooth projective surfaces of Kodaira dimension zero.
\end{abstract}

\maketitle

\blfootnote{\scriptsize{Keywords: elliptic surfaces, positive characteristic, Picard schemes, finite group schemes, Artin-Schreier covers}}

\blfootnote{\scriptsize{2020 Mathematics Subject Classification: Primary: 14J27, 14K30, Secondary: 14L15, 14L30.}}

\setcounter{tocdepth}{1} 
\tableofcontents

\section{Introduction}

The present article is devoted to the study of the Picard scheme of certain isotrivial elliptic surfaces; in particular, to determine whether they are smooth or singular. By Cartier's theorem, they may be singular only in positive characteristic.

\medskip

First, let us recall Igusa's construction, which works over any algebraically closed field $k$ of characteristic $p=2$. We fix an ordinary elliptic curve $E$, so that the two-torsion subgroup $E[2]\subset E$ contains a copy of $\mathbb{Z}/2\mathbb{Z}$, which also acts on any other elliptic curve $X$ by the involution $\iota\colon x\mapsto -x$, and the quotient $X/(\mathbb{Z}/2\mathbb{Z})$ is isomorphic to $\P^1$.
This defines a diagonal action of $\mathbb{Z}/2\mathbb{Z}$ on the product $E\times X$, acting by translations on $E$ and by the involution $\iota$ on $X$. The induced quotient 
\[
E\times^{(\mathbb{Z}/2\mathbb{Z})} X,
\]
is called an \emph{Igusa surface}. Introduced in \cite{Igusa}, these were the first examples of smooth projective surfaces having a singular Picard scheme. Notice that Igusa surfaces are bielliptic surfaces; in particular, they have Kodaira dimension zero. Note also that \(\mathbb{Z}/2\mathbb{Z}\) is unipotent in characteristic \(2\), and that the action on \(X\) has fixed points, namely \(X[2]\).

\medskip

Later, various examples of smooth projective surfaces with singular Picard schemes were constructed for arbitrary characteristic $p>0$.
As a direct generalization of Igusa surfaces, Katsura and Ueno considered the analogous construction in which the double cover $X\to \P^1$ is replaced with an Artin-Schreier cover of $\P^1$, and more generally $(\mathbb{Z}/p\mathbb{Z})^n$-covers over $\P^1$ (see \cite[\S 8]{Katsura_Ueno}). 
Another class of examples is obtained by considering quotients by finite group schemes, which may themselves be non-reduced. In Kodaira dimension zero, this yields bielliptic and quasi-bielliptic surfaces, see \cite{BMIII,BMII}. Both constructions produce elliptic surfaces with singular Picard schemes, as pointed out in \cite[Theorem 1.2, Proposition 2.3]{Liedtke_note}. 

\medskip

Observe that, in both of the above generalizations, the resulting smooth projective surfaces have Kodaira dimension at most one. There also exist examples of smooth projective surfaces of general type with singular Picard schemes. Such examples were found by Serre, see \cite[Propositions 15 and 16]{Serre} and \cite[Proposition 3.1]{Liedtke_note}.

\subsection{Isotrivial elliptic surfaces via equivariantly normal curves}

From now on, we fix an algebraically closed field $k$ of arbitrary characteristic $p>0$. All objects and morphisms, unless otherwise specified, are assumed to be defined over $k$. Group schemes are assumed to be of finite type.

\medskip

Our main object of study is the Picard scheme of an isotrivial elliptic surface $S$, which is isomorphic to a contracted product
\begin{align}\tag{$\star$}
    \label{eq:contracted product}
E\times^G X := (E\times X)/G,
\end{align}
where $E$ is an elliptic curve, $G\subset E$ is a finite subgroup scheme, and $X$ is a projective $G$-normal curve; i.e., every finite birational $G$-equivariant morphism with target $X$ is an isomorphism. These smooth projective surfaces have Kodaira dimension one if the arithmetic genus of \(X\) is at least two (see \cite[Theorem A]{fong_maccan}). The theory of equivariantly normal curves (or $G$-normal curves) was introduced by Brion in \cite{G-normal-curves}. It provides an elegant remedy to the fact that the action of $G$, when $G$ is non-reduced, might not lift to the normalization of the variety on which $G$ acts.

\medskip

For our purposes, this provides a unified framework for studying smooth projective surfaces with singular Picard scheme. Indeed, in \Cref{prop:Rosenlicht_application}, we show that every smooth projective surface $S$ with a singular $\Pic_S$ and a positive-dimensional automorphism group is necessarily isomorphic to a contracted product as above. We call this a \emph{generalized Igusa surface} (see \Cref{def:generalized_Igusa}). In particular, this construction recovers the examples of Igusa, Katsura-Ueno, and Bombieri-Mumford. 

\subsection{The main results}

To determine whether the Picard scheme of a smooth projective surface $S$ is smooth, a natural approach is to compute the birational invariant
\[
\delta_I(S) = h^1(S,\O_S) - \frac{b_1(S)}{2} \geq 0,
\]
which vanishes if and only if $\Pic_S$ is smooth (it is the difference of the dimension of the Lie algebra of the Picard scheme and the dimension of the Picard scheme).

\medskip

Following \cite{BMII}, for an elliptic fibration $f\colon S\to Y$, we have
\[
R^1f_*(\O_S) = \mathcal{L} \oplus \mathcal{T},
\]
where $\mathcal{L}$ is an invertible sheaf and $\mathcal{T}$ is a torsion sheaf. A fiber $F$ over a point in the support of $\mathcal{T}$ is called a \emph{wild} fiber (or equivalently, $h^0(F,\O_F) \geq 2$). 

\medskip

Combining the theories of equivariantly normal curves and elliptic surfaces, we obtain the following result:

\begin{theoremA}\label{thm:A}
    Let $G$ be a finite subgroup scheme of an elliptic curve $E$ acting faithfully on a projective $G$-normal curve  $X$, with geometric quotient $\pi\colon X\to Y$. Let $S=E\times^G X$, which is a smooth projective surface, equipped with the isotrivial elliptic fibration
    $f\colon S \to Y$.
    Then 
    \[
    \delta_I(S)=h^0(X,\omega_X)^G - g(Y),
    \]
    where $\omega_X$ denotes the dualizing sheaf of $X$.
    Moreover, the following hold:
    \begin{enumerate}
        \item The multiplicity of a fiber $f^{-1}(\pi(x))$ is equal to the order of $\mathrm{Stab}_G(x)$.
        \item If $x\in X$ has a nontrivial unipotent stabilizer, then $f^{-1}(\pi(x))$ is a wild fiber.
        \item\label{thm:A.iii} $S$ is a generalized Igusa surface if and only if $h^0(Y,\mathcal{T})\geq 2$. In particular, these two equivalent conditions hold if \(f:S\to Y\) admits two or more wild fibers, or if it has a unique wild fiber \(F\) such that \(h^0(F,\mathcal{O}_F)>2\).
    \end{enumerate}
\end{theoremA}

\Cref{thm:A} shows that for our smooth projective surfaces $S=E\times^G X$, the smoothness of $\Pic_S$ is entirely determined by the datum of the $G$-normal curve $X$. Moreover, it is known that any elliptic surface with at least two wild fibers has a singular Picard scheme (see e.g. \cite[Proposition 2.1]{Liedtke_note}). \Cref{thm:A} (\ref{thm:A.iii}) gives a refinement of this statement for our surfaces. As a recipe for constructing isotrivial elliptic surfaces with singular Picard schemes, we get:

\begin{corollaryA}\label{coro:B}
    Assume moreover that $G$ is unipotent in \Cref{thm:A}. Then $S$ is a generalized Igusa surface, i.e., $\Pic_S$ is singular, if one of the following conditions hold:
    \begin{enumerate}
        \item There exist $x_1,x_2\in X$ with distinct orbits and nontrivial stabilizers. 
        \item The arithmetic genus of \(X\) is greater than or equal to one and \(Y\cong \mathbb{P}^1\). In particular, the action of \(G\) on \(X\) is never free in this case.
    \end{enumerate}
\end{corollaryA}

Taking $G$ to be diagonalizable never produces an isotrivial elliptic surface $S$ with a singular Picard scheme (see \Cref{thm:diag_reduced_pic}). In contrast, choosing \(G\) unipotent, and $p_a(X)\geq 1$ such that the quotient \(Y=X/G\) is isomorphic to \(\mathbb{P}^1\), then $\Pic_S$ is automatically singular. In broad terms, the two opposite classes of finite group schemes, namely diagonalizable and unipotent, lead to opposite behaviors of the Picard schemes of the resulting isotrivial elliptic surfaces.

\medskip

Next, we give applications of \Cref{thm:A} and \Cref{coro:B}. Since $G$ is a subgroup scheme of an elliptic curve $E$, there are two natural choices for unipotent $G$, which are minimal for the inclusion: the constant group scheme $\mathbb{Z}/p\mathbb{Z}$ and the infinitesimal group scheme $\aalpha_p$, corresponding respectively to the cases where $E$ is ordinary and supersingular.

\medskip

For $E$ ordinary and $G=\mathbb{Z}/p\mathbb{Z}$, the $G$-normal curve $X$ is smooth, and the quotient morphism $\pi\colon X\to Y$ is an Artin-Schreier cover. The induced contracted product is a direct generalization of the examples of Igusa and Katsura-Ueno, where $Y=\P^1$. Combining \Cref{thm:A} with a result of Valentini and Madan (see \cite{Valentini_Madan}), which allows us to compute the $G$-invariants of differential forms, we obtain the following statement:

\begin{theoremA}\label{thm:C} 
    Let $E$ be an ordinary elliptic curve, let $\pi\colon X\to Y$ be an Artin-Schreier cover, and let $S=E\times^{\mathbb{Z}/p\mathbb{Z}} X$ be the induced contracted product, which is equipped with the elliptic fibration $f\colon S \to Y$. If $\pi$ is unramified then $\Pic_S$ is smooth. 
    Otherwise, we denote by $\{q_i\}_{i=1,\ldots,s}$ the branched locus, with each $q_i$ of ramification jump $m_i>0$. Then each fiber $f^{-1}(q_i)$ is wild, and 
    \[
    \delta_I(S) =   \left( \sum_{i=1}^s \left\lfloor \frac{(m_i+1)(p-1)}{p} \right\rfloor \right) -1,
    \]
    and $\Pic_S$ is smooth if and only if \(\pi\colon X\to Y\) has at most a single branch point with ramification jump one.
\end{theoremA}

In \Cref{thm:C}, notice that for branched point $q_i$, we have that 
\[
\left\lfloor \frac{(m_i+1)(p-1)}{p} \right\rfloor \geq 1.
\]
Note that the explicit formula for \(\delta_I(S)\) for a generalized Igusa surface associated to an Artin-Schreier cover gives a new class of examples of elliptic fibrations with arbitrarily large \(\delta_I(S)\), giving a different proof of \cite[Proposition 2.3]{Liedtke_note}. Moreover, our construction can also be used to produce elliptic fibrations with large \(\delta_I(S)\) over an arbitrary smooth projective curve, not just over \(\mathbb{P}^1\).

\medskip

As another application, we consider the case where $E$ is supersingular and $G=\aalpha_p$. To the extent of the authors' knowledge, there is no classification of $\aalpha_p$-normal curves, and one can expect such a classification problem to be difficult (see e.g. \cite{ClassificationMuPNormal} for a classification of $\muu_p$-normal curves). However, 
interesting examples of $\aalpha_p$-normal curves have recently been constructed by Hilario and Schröer in \cite{Generalized-quasielliptic}. They introduce a generalization of the rational cuspidal plane curve $y^2 = x^3$, which we denote by $X_p$ (see \Cref{subsection:generalized_Igusa_alpha_p}). Analogously to the case where $G=\mathbb{Z}/p\mathbb{Z}$ is constant, we consider the induced isotrivial elliptic surface $S_p=E\times^{\aalpha_p} X_p$, and we obtain:

\begin{propositionA}\label{prop:E}
    Let $E$ be a supersingular elliptic curve, let 
    $X_p$ be the rational curve defined by Hilario-Schröer, and let $S_p=E\times^{\aalpha_p} X_p$ be the induced contracted product. Then 
    \[
    \delta_I(S_p) = p-2,
    \]
    and $S_p$ is a generalized Igusa surface if and only if $p\geq 3$.
    Moreover, the following hold:
    \begin{enumerate}
        \item If $p=3$, then $S_p$ is a quasi-bielliptic surface.
        \item If $p\geq 5$, then $\kappa(S_p)=1$.
    \end{enumerate}
\end{propositionA}

It remains to describe the group scheme structure of $\Pic_S^0$, where $S$ is a generalized Igusa surface. 
For most of the cases where \(S\) has Kodaira dimension zero, this follows from the PhD thesis of Kroon \cite[Table 7.1]{Kroon}. Nonetheless, a few cases of their classification are left conjectural. The below result is a partial answer to these conjectures, see the discussion following \Cref{prop: picinf quasi bielliptic}.

Our strategy is based on the results of Jensen \cite{Jensen} (see also \Cref{thm:Jensen} and \Cref{cor:Jensen}). Under some additional assumptions on the group $G$, this amounts to describing the $G$-invariant closed subgroup scheme $(\Pic_X^0)^G \subset \Pic_X$, where $X$ is the $G$-normal curve. This way, we get a full description of \((\Pic^0_S)_{\mathrm{inf}}\), namely the quotient of the Picard scheme by its reduced subscheme, in Kodaira dimension zero as follows (see \Cref{prop: picinf quasi bielliptic,prop: picinf bielliptic}).

\begin{theoremA}\label{prop:F}
    Let \(S=E\times^GX\) be a generalized Igusa surface of Kodaira dimension \(0\). Then \((\Pic^0_S)_{\mathrm{inf}}\) is isomorphic to one of the following group schemes:
    \begin{enumerate}
        \item \(\aalpha_3\), where \(X\) is either a rational curve with a cusp, or a supersingular elliptic curve (\(p=3\));
        \item \(\aalpha_2\), where \(X\) is either a rational curve with a cusp, or a supersingular elliptic curve (\(p=2\));
        \item  \(\aalpha_4\), where $X$ is a rational curve with a cusp (\(p=2\));
        \item $\muu_{2}$, where $X$ is an ordinary elliptic curve (\(p=2\));
        \item $X[2]$, where $X$ is a supersingular elliptic curve  (\(p=2\)). 
    \end{enumerate}
\end{theoremA}

By \cite[Theorems 1.1 and 1.2]{Liedtke_note}, the smooth projective surfaces with Kodaira dimension $\kappa(S)\leq 0$ and having a singular Picard scheme  are either (quasi)-bielliptic surfaces or non-classical Enriques surfaces. The latter only occurs in characteristic two, and their Picard schemes are completely described in \cite[Theorem 1.2.1]{Cossec-Dolgachev-Liedtke}. So the above theorem completes the description of the non-reduced structure of singular Picard schemes for all surfaces of Kodaira dimension \(\leq 0\).

\medskip

In Kodaira dimension one, generalized Igusa surfaces are more abundant, see e.g \Cref{thm:C} and \Cref{prop:E}. At this stage, we do not see how to describe \((\Pic^0_S)_{\mathrm{inf}}\) uniformly. We do however provide a description of the natural \(\aalpha_p\)-action on \(\Pic_{X_p}\) in \Cref{lem: picardXp}, and show that computing \((\Pic_{S_p})_{\mathrm{inf}}\) reduces to computing the fixed point scheme for this action (see \Cref{cor:Jensen}).

\subsection{Structure of the article}

The article is organized as follows. 

\medskip

In \Cref{section:preliminaries}, we recall the key notion of $G$-normal curves, where $G$ is a finite group scheme, as well as their properties. We also introduce here one of the central objects in this article: the isotrivial elliptic surfaces $S = E\times^G X,$ where $E$ is an elliptic curve, $G\subset E$ is a finite subgroup scheme, and $X$ is a projective $G$-normal curve.
We show that this family of smooth projective surfaces contains the ones with singular Picard schemes and positive-dimensional automorphism groups (\Cref{prop:Rosenlicht_application}).
We also briefly recall the construction of the Picard scheme and the canonical bundle formula for elliptic surfaces. In the end of \Cref{section:preliminaries}, we fix the notations used in the rest of the paper.

\medskip

In \Cref{section:Picard_scheme}, we revisit several well-known formulas for elliptic surfaces in our setting, which makes them simpler and more explicit. Using the theory of $G$-normal curves, we get the following new formula (\Cref{lem:invariants}): 
\[
\delta_I(S) = h^0(\omega_X)^G - g(Y).
\]
We also precise the relationship between the reducedness of $\Pic_S$ and the existence of wild fibers (\Cref{prop:one_wild_fiber}). A useful result is \Cref{prop:stab_multiplicity}, where we show that the multiple fibers of $f\colon S\to Y$ correspond to the images under $\pi\colon X\to Y$ of points with nontrivial $G$-stabilizers. In the second half of \Cref{section:Picard_scheme}, we focus on the case where $G$ is unipotent. We prove that points of $X$ with unipotent stabilizers necessarily give rise to wild fibers. Combining the results established in this section, we prove \Cref{thm:A} and \Cref{coro:B}.

\medskip

\Cref{section:applications} discusses some applications. 
The case $G=\mathbb{Z}/p\mathbb{Z}$ shows that any Artin-Schreier cover with enough ramification has an associated surface with singular Picard scheme (\Cref{thm:C}),
while the case $G=\boldsymbol{\alpha}_p$ leads to \Cref{prop:E}. 

\medskip

In \Cref{sec: picinf}, we describe the group scheme structure of the singular Picard schemes of (quasi)-bielliptic surfaces. More precisely, we determine the quotient of the Picard group by its reduced subscheme, as an infinitesimal (topologically trivial) group scheme. 

\medskip

\noindent \textbf{Acknowledgments.}
We are grateful to Fabio Bernasconi, Michel Brion, Raymond Cheng, Bianca Gouthier, Amine Koubaa, and Gebhard Martin for helpful discussions. We would like to particularly thank Bianca Gouthier for her contributions to this project and for her involvement in its earlier stages.
The second author acknowledges the support of the DFG, through the research grants Le 3093/5-1 and Le 3093/7-1 (project numbers 530132094; 550535392). 

\medskip

\noindent \textbf{AI disclosure.} AI-assisted tools (ChatGPT and Claude) were used for proofreading and for identifying potentially relevant references (in particular \cite{Kroon} and \cite{Valentini_Madan}). All statements and proofs were written by the authors, who take full responsibility for the content of this article.

\section{Preliminaries}\label{section:preliminaries}

\subsection{Isotrivial elliptic surfaces via equivariantly normal curves}

First, let us recall the notion of equivariantly normal curves, introduced by Brion in \cite{G-normal-curves}.

\begin{definition}\cite[Definition 4.1]{G-normal-curves}\label{def:G-normal}
    Let $G$ be a finite group scheme. A curve $X$ with an action of $G$ is $G$-normal if every finite birational $G$-equivariant morphism with target $X$ is an isomorphism.
\end{definition}

\begin{lemma}\label{lem:faithful->gen.free}
    Let \(G\subset E\) be a finite commutative group scheme with Lie algebra of dimension \(\leq 1\) (e.g. a finite subgroup scheme of an elliptic curve). If \(G\) acts faithfully on a curve \(X\), then it acts generically freely.
\end{lemma}

\begin{proof}
We have \(G\cong G_{\mathrm{red}}\times G/G_{\mathrm{red}}\) by \cite[Theorem 15.5]{FiniteGruopSchemesPink}, where \(G_{\mathrm{red}}\) is a finite constant group acting faithfully on the curve \(X\), thus it acts generically freely.
    Since \(\dim_k\mathrm{Lie}(G)\) has dimension at most \(1\), it follows that it is trivial or acts faithfully on \(X\) by \cite[Lemma 5.3]{G-normal-curves}.
\end{proof}

\begin{proposition}[{\cite[Corollary 4.16]{G-normal-curves}}]\label{prop:quotient_G-normal}
    Let $G$ be a finite group scheme and let $X$ be a $G$-normal curve. Then the quotient $Y=X/G$ is a smooth curve.
\end{proposition}

As a direct consequence of the above proposition, we have:

\begin{corollary}\label{cor:fin-flat}
    Let \(X\) be a \(G\)-normal curve, then the quotient morphism \(\pi:X\to Y\) is finite flat. If \(\mathrm{Lie}(G)\) is one-dimensional and the action is faithful, then the degree of \(\pi\) is \(|G|\).
\end{corollary}
\begin{proof}
    Since \(\pi:X\to Y=X/G\) is a geometric quotient, it is dominant, and \(Y\) is smooth by the above \Cref{prop:quotient_G-normal}. Hence, the statement follows from \cite[\S 4.3.1. Proposition 3.9]{Liu}.
    To see the last statement, note that by \Cref{lem:faithful->gen.free}, the map \(X\to Y\) is generically a \(G\)-torsor, hence its generic degree is \(|G|\), and thus \(\pi\) has degree \(|G|\) by flatness.
\end{proof}

In the following proposition, we introduce the \emph{contracted product} 
\[
E\times^G X.
\]
By this notation, we mean the quotient of $E\times X$ by a diagonal action of $G$: for every $g\in G$ and $(p,x)\in E\times X$, it is given by 
\[
(p,x) \mapsto (p\cdot g^{-1}, g\cdot x).
\]

\begin{proposition}[{\cite[Proposition 5.6]{G-normal-curves}}] \label{prop:isotrivial_surfaces}
    Let $S$ be a smooth projective surface equipped with a faithful action of an elliptic curve $E$. Then there exists a finite subgroup scheme $G\subset E$ acting faithfully on a projective $G$-normal curve $X$ such that $S = E\times^G X$. 
\end{proposition}

\begin{remark}
In \Cref{def:G-normal}, the $G$-action on $X$ is not assumed to be faithful in general. However, for the purposes of the present text, we can assume this action being faithful by \Cref{prop:isotrivial_surfaces}, as our $G$-normal curve $X$ arises from a smooth projective surface $S$ with a faithful action of an elliptic curve.

The faithfulness of the action on $X$ is not part of the statement of \cite[Proposition 5.6]{G-normal-curves}, but it is part of the proof. For self-containedness, let us explain this short argument. Writing $S=E\times^G X$, the finite group subscheme $G$ acts by translations on the first factor of the contracted product as 
\[
g\cdot (p,x) = (g\cdot p,x).
\]
Using that $E$ is commutative and the equivalence relation defining the contracted product, we obtain
\[
 g\cdot (p,x)= (p,g\cdot x).
\]
Since $G$ acts faithfully on $S$, it follows from the above equality that $G$ acts faithfully on $X$ as well.
\end{remark}

\subsection{The Picard scheme} 

For a variety $Z$, we denote as $\mathrm{Pic}(Z)$ the group of the line bundles on $Z$, up to linear equivalence. The \emph{Picard scheme} of $Z$ is (when defined) the scheme $\underline{\mathrm{Pic}}_Z$ representing the relative Picard functor (over $k$) given by
\[
T \leadsto \mathrm{Pic}(Z \times T)/\mathrm{Pic}(T),
\]
where $T$ varies among the $k$-schemes. If $Z$ is a projective variety, then its Picard scheme exists, and it is a group scheme, locally of finite type over $k$ (not necessarily reduced). Its connected component of the identity $\underline{\mathrm{Pic}}^0_Z$ is an algebraic group. For these facts and more on the definition of the different Picard functors, see \cite[Section 4 and 5]{Kleiman}.

\medskip

In fact, the surfaces $S=E\times^G X$ are a key source of examples of smooth projective surfaces with singular Picard schemes. We have the following proposition:

\begin{proposition}\label{prop:Rosenlicht_application}
    Let $S$ be a smooth projective surface whose automorphism group scheme $\Aut_S$ has positive dimension and whose Picard scheme $\Pic_S$ is singular. Then $E=(\Aut^0_S)_{\mathrm{red}}$ is an elliptic curve, and there exists a finite subgroup scheme $G\subset E$ acting faithfully on a projective $G$-normal curve $X$ such that 
        \[
        S = E\times^G X.
        \]
\end{proposition}

\begin{proof}
    By Chevalley's structure theorem, there exists a short exact sequence
    \[
    1 \to L \to (\Aut^0_S)_{\mathrm{red}} \to A \to 1,
    \]
    where $L$ is a connected linear algebraic group and $A$ is an abelian variety. If $L$ is not trivial, it follows that $S$ is birationally equivalent to $C\times \P^1$, where $C$ is a smooth projective curve (see \cite[Theorem 10]{Rosenlicht}, and also \cite[Proposition 2.5.1]{BFT}). Then $h^1(S,\O_S) = b_1(S)/2 = g(C)$,
    and this contradicts $\Pic_S$ being singular. Therefore, $(\Aut^0_S)_{\mathrm{red}}=A$ is an abelian variety. 

    Assume that there exists an open orbit $A\cdot x$, then the isomorphism
    \[
    A\cdot x \cong A /\mathrm{Stab}(x) 
    \]
    implies that $A$ has dimension two, as $\mathrm{Stab}(x) $ is finite (see e.g. \cite[Proposition 2.2.1]{BSU}). Then $S$ is an abelian surface acting on itself by translations, and it also follows that $\Pic_S$ is reduced (see e.g. \cite[Theorem 1.2]{Liedtke_note}). This leads again to a contradiction. Since $\Aut_S$ has positive dimension, it remains to consider the case where $(\Aut^0_S)_{\mathrm{red}}$ is isomorphic to an elliptic curve $E$ and we conclude with \Cref{prop:isotrivial_surfaces}.
\end{proof}

This motivates the following definition:

\begin{definition}\label{def:generalized_Igusa}
    A \emph{generalized Igusa surface} is a smooth projective surface with a positive dimensional automorphism group and a singular Picard scheme.
\end{definition}

In previous work by two of the authors, a preliminary study on Picard schemes of contracted products \(S=E\times^GX\) was initiated, showing that if \(G\) is diagonalizable then \(S\) is not a generalized Igusa surface:

\begin{theorem}[{\cite[Theorem B]{fong_maccan}}]\label{thm:diag_reduced_pic}
    Let $S=E\times^G X$, where $E$ is an elliptic curve and $G$ is a finite \emph{diagonalizable} subgroup scheme of $E$ acting on the projective $G$-normal curve $X$. Then the Picard scheme $\Pic_S$ is reduced.
\end{theorem}

\subsection{The canonical bundle formula} 

The contracted product $S=E\times^G X$ is equipped with a morphism $f\colon S\to Y=X/G$, which is an elliptic fibration admitting finitely many \emph{multiple fibers}, i.e., non-reduced schematic fibers. The following definition and theorem are classical, see e.g. \cite{BMIII,BMII} for a standard reference:

\begin{definition}
    Let $f\colon S\to Y$ be an elliptic fibration as above and let $F\subset S$ be a fiber. 
    \begin{enumerate}
        \item The multiplicity of $F$ is the integer $m\geq 1$ such that 
        \[
            F = m\cdot F_{\mathrm{red}},
        \]
        as Cartier divisor, and where $F_{\mathrm{red}}\subset F$ denotes the underlying reduced subscheme. In particular, $m\geq 2$ if and only if $F$ is a multiple fiber.
        \item A multiple fiber $F$ is a \emph{tame} fiber if $h^0(F,\O_F)=1$. 
        \item A multiple fiber $F$ is a \emph{wild fiber} if $h^0(F,\O_F)\geq 2$.
    \end{enumerate}
\end{definition}

Let us recall the canonical bundle formula, which is a fundamental result in the study of elliptic surfaces. Let $f\colon S\to Y$ be an elliptic fibration. Following \cite{BMII}, we write 
\[R^1f_*(\O_S) = \mathcal{L} \oplus \mathcal{T},\]
where $\mathcal{L}$ is a line bundle and $\mathcal{T}$ is supported in finitely many points. The fibers above the support of $\mathcal{T}$ are precisely the wild fibers.

\begin{theorem}[{\cite[Theorem 2]{BMII}}]\label{thm:canonical_bundle_formula}
    With the notations above,
    \[
    \omega_S = f^*(\mathcal{L}^{-1} \otimes \omega_Y) \otimes \O_S\left(\sum_\lambda a_\lambda P_\lambda\right),
    \]
    where 
    \begin{enumerate}
        \item $m_\lambda P_\lambda$ are the multiple fibers,
        \item $0 \leq a_\lambda < m_\lambda$,
        \item $a_\lambda = m_\lambda - 1$ if $m_\lambda P_\lambda$ is not wild,
        \item $\deg(\mathcal{L}^{-1}) = \chi(\mathcal{O}_S)+ h^0(Y,\mathcal{T})$.
    \end{enumerate}
\end{theorem}

The following proposition is well-known to experts, but we could not locate any reference. For self-containedness, we prove the statement:

\begin{proposition}\label{prop:sections_wild}
    Let $f\colon S\to Y$ be an elliptic surface over a smooth projective curve $Y$. Then
    \[
        h^0(Y,\mathcal{T}) \geq \sum_F \big(h^0(F,\O_F) -1\big),
    \]
    where the sum runs over the finite set of wild fibers.
\end{proposition}

\begin{proof}
    Since $\mathcal{T}$ is a torsion sheaf, we can write 
    \[
    \mathcal{T} = \bigoplus_{y} \mathcal{T}^{(y)},
    \]
    where the direct sum runs over the support of the wild fibers, which is a finite set, and $\mathcal{T}^{(y)}$ is a skyscraper sheaf supported at $y\in Y$.
    
    Consider the wild fiber $F$ at $y\in Y$ and the following cartesian square
    \[
    \begin{tikzcd}
	F \ar[r,hook]\ar[d,"f_y"] & S \ar[d,"f"] \\
	{\Spec(k)} \ar[r,hook,"y"]& Y.
	\end{tikzcd}
    \]
    Since all fibers of \(f\) are one-dimensional and \(f\) is proper, by proper base change in top degree (see the first two lemmas of \cite[\href{https://stacks.math.columbia.edu/tag/0EX2}{Section 0EX2}]{stacks-project}), we have that 
    \[
    R^1f_*(\O_S) \otimes_{\O_Y} k(y) \cong H^1(F,\O_F) .
    \]
    The left hand side equals $H^0(\Spec k,\mathcal{L}|_{y} \oplus \mathcal{T}|_{y}) \cong k \oplus (\mathcal{T}_y/\mathfrak{m}_y\mathcal{T}_y)$, where $\mathfrak{m}_y$ denotes the unique maximal ideal of $\O_{Y,y}$. In particular, this has dimension at most $1+ h^0(Y,\mathcal{T}^{(y)})$ over \(k\), with equality if and only if \(\mathcal{T}^{(y)}\cong y_*k^{\oplus n}\) for some positive integer \(n\). 
    
    Moreover, $h^1(F,\O_F)=h^0(F,\O_F)$. Indeed, we have that
    \[
    \chi (\O_F)= \chi(\O_S) - \chi(\O_S(-F)) = -\frac{1}{2}F\cdot(F+K_S) =0.
    \]
    The first equality follows from the short exact sequence $0\to \O_S(-F) \to \O_S \to \O_F \to 0$ and the additivity of the Euler characteristic, the middle one is Riemann-Roch formula, while the last equality follows from $F^2=0$ and the fact that $F\cdot K_S=0$ (see \Cref{thm:canonical_bundle_formula}). 
    
    This proves that $h^0(F,\O_F) -1 \leq h^0(Y,\mathcal{T}^{(y)})$. Summing the contribution of all wild fibers, we obtain the stated inequality.
\end{proof}

\begin{remark}
    For every fiber $F$ of the elliptic surface $f\colon S\to Y$, we have the inequalities
    \[
    1 \leq h^0(F,\O_F) \leq m,
    \]
    where $m\geq 1$ denotes the multiplicity of the fiber $F$. Indeed, $F=mF_{\mathrm{red}}$ as Cartier divisors, and taking the global sections in the short exact sequence
    \[
    0 \to \O_{F_{\mathrm{red}}}(-(m-1)F_{\mathrm{red}}) \to \O_{mF_{\mathrm{red}}} \to \O_{(m-1)F_{\mathrm{red}}} \to 0,
    \]
    we get that $h^0(mF_{\mathrm{red}},\O_{mF_{\mathrm{red}}}) \leq 1+h^0((m-1)F_{\mathrm{red}},\O_{(m-1)F_{\mathrm{red}}})$ as $F_{red}^2=0$. Then by induction, we see that $m$ is an upper bound of $h^0(F,\O_F)$.
\end{remark}

\subsection{Setup}\label{sub:setup} Let us fix the notations used throughout this article. We denote by $E$ an elliptic curve, $G\subset E$ a finite subgroup scheme acting faithfully on a projective $G$-normal curve $X$. By \cite[Proposition 4.17]{G-normal-curves}, the contracted product 
$$ S=E\times^G X $$ is a smooth projective surface, equipped with two projections: one onto the elliptic curve $E/G$ which is isogeneous to $E$, and the second one onto the smooth projective curve $Y=X/G$. This gives the following commutative diagram:
\[
\begin{tikzcd}
    E \times X \arrow[dd, "pr_X"] \arrow[r,"q"] & S\arrow[dd,"f"] \arrow[r,"h"] & E/G\\
    &&\\
    X \arrow[r, "\pi"] & Y &.
\end{tikzcd}
\]
The quotient morphism $q$ is a $G-$torsor. The map $h$ is a surjective morphism locally trivial for the $fppf$ topology and with fiber $X$. On the other hand, the morphism $f$ is an elliptic fibration which can admit \emph{multiple fibers}, i.e., schematic fibers that are not reduced.

\section{Picard scheme of isotrivial elliptic surfaces}\label{section:Picard_scheme}

\subsection{The birational invariant $\delta_I(S)$}

Two cohomological invariants play a central role in the smoothness of the Picard scheme: $h^1(S,\O_S)$ and $b_1(S)$. The first one, called \emph{the irregularity}, is the dimension of $\mathrm{Lie}(\Pic_S)$; the second one, which is the first Betti number, is the double of the dimension of $\Pic_S$. In particular, the non-negative integer
\[
\delta_I(S) = h^1(S,\O_S) - b_1(S)/2 \geq 0
\]
measures the obstruction for the group scheme $\Pic_S$ to be reduced, and it is zero if and only if $\Pic_S$ is reduced. Notice that $\delta_I(S)$ is a birational invariant for smooth projective surfaces, as blowing up a smooth point does not affect $\Pic^0_S$. Let us recall that by \cite[Theorem A]{fong_maccan}, we already have
\[
\begin{cases}
    b_1(S) = 2 + 2g(Y) \\ b_2(S) = 2 + 4g(Y)
\end{cases}
\]
where $Y$ denotes the smooth projective curve $X/G$.

\begin{lemma}\label{lem:invariants}
    For the surfaces $S = E\times^G X$, the following equalities hold:
        \begin{enumerate}
            \item\label{lem:invariants.1} $\chi(\O_S)=0$,
            \item\label{lem:invariants.2} $h^1(S,\O_S) = 1 + h^0(X,\omega_X)^G$,
            \item\label{lem:invariants.3} $\delta_I(S) = h^0(X,\omega_X)^G - g(Y) = h^0(Y,\mathcal{L})+h^0(Y,\mathcal{T}) -1$.
        \end{enumerate}
    Moreover, we have $\deg(\mathcal{L})=-h^0(Y,\mathcal{T})\leq0$. The elliptic fibration $f\colon S\to Y$ has no wild fibers if and only if $\mathcal{L} =\O_Y$; otherwise, $\deg(\mathcal{L})<0$. 
\end{lemma}

\begin{proof}
    By the canonical bundle formula, we have that $K_S^2 = 0$. Using that $b_1(S) = 2 + 2g(Y)$ and $b_2(S) = 2 + 4g(Y)$ (see \cite[Theorem A]{fong_maccan}), we get that the topological Euler characteristic of $S$ vanishes. Then (1) follows from Noether's formula. 
    
    Using (1) and Serre duality, we get $h^1(S,\O_S) =  1 + p_g(S)$, and $p_g(S) = h^0(X,\omega_X)^G$ by \cite[Lemma 3.25]{fong_maccan}. Together, this yields (2).
    
    Next we prove (3). Subtracting both sides of (2) by $b_1(S)/2$, that equals $1 +g(Y)$, we get that $\delta_I(S)= h^0(X,\omega_X)^G - g(Y)$. As in \cite{BMII}, using the Grothendieck-Leray spectral sequence, one derives the short exact sequence
    \[
        0 \to H^1(Y,\O_Y) \to H^1(S,\O_S) \to H^0(Y,R^1f_*(\O_S)) \to 0.
    \]
    Therefore, $h^0(Y,R^1f_*(\O_S)) = h^1(S,\O_S) - h^1(Y,\O_Y)$. Since $R^1f_*(\O_S) = \mathcal{L} \oplus \mathcal{T}$ and $h^1(Y,\O_Y)=b_1(S)/2-1$, we obtain the last equality of (3). 
    
    We now prove the last assertions. By \cite[Theorem 2]{BMII}, we have $\deg(\mathcal{L}^{-1}) = \chi(\O_S) + h^0(Y,\mathcal{T})$. Since $\chi(\O_S)=0$, it follows that $\deg(\mathcal{L}) = -h^0(Y,\mathcal{T}) \leq 0$. In particular, $\deg(\mathcal{L})=0$ if and only if $f$ has no wild fibers. In this case, $h^0(Y,\mathcal{T})=0$; hence, $\delta_I(S) = h^0(Y,\mathcal{L}) - 1$. As $\delta_I(S)\geq 0$, we obtain that $h^0(Y,\mathcal{L})\geq 1$, so $\mathcal{L}$ is trivial. 
\end{proof}

In particular, as proven in \cite[Proposition 2.1 (2)]{Liedtke_note}, we see from \Cref{lem:invariants} that if $f\colon S\to Y$ has $w\geq 2$ wild fibers, then $\delta_I(S) \geq w-1$, and $S$ is a generalized Igusa surface. For our surfaces $S=E\times^G X$, we obtain an analogue of \cite[Proposition 2.1 (1)]{Liedtke_note}: $\Pic_S$ is reduced if there is no wild fibers. Moreover, if there exists a unique wild fiber, we show that $\Pic_S$ is singular if this fiber is "wild enough".

\begin{proposition}\label{prop:one_wild_fiber}  
    The surface $S = E\times^G X$ is a generalized Igusa surface if and only if 
    \[
    h^0(Y,\mathcal{T})\geq 2.
    \]
    Let $w\geq 0$ be the number of wild fibers of $f\colon S\to Y$. Then the following hold:
    \begin{enumerate}
        \item If $w=0$, then $\Pic_S$ is reduced.
        \item If $w=1$, then 
        \[
        \delta_I(S)\geq  h^0(F,\O_F) -2,
        \]
        where $F$ denotes the unique wild fiber. In particular, $S$ is a generalized Igusa surface if $h^0(F,\O_F)>2$. 
        \item If $w\geq 2$, then $S$ is a generalized Igusa surface.
    \end{enumerate}
\end{proposition}

\begin{proof} 
By \Cref{lem:invariants}, we have $\delta_I(S) = h^0(Y,\mathcal{L})+h^0(Y,\mathcal{T}) -1$, which is positive if $h^0(Y,\mathcal{T})\geq 2$. Conversely, if $\delta_I(S) > 0$, then $f$ admits at least one wild fiber and $h^0(Y,\mathcal{L})=0$. This implies that $h^0(Y,\mathcal{T}) \geq 2$. 

(1) If $f$ has no wild fiber, then $h^0(Y,\mathcal{T}) =0$ and $\mathcal{L}=\mathcal{O}_Y$ by \Cref{lem:invariants}. In particular, $\delta_I(S)=0$, i.e., $\underline{\mathrm{Pic}}_{S}$ is reduced.

(2) Assume that $w=1$. Then $h^0(Y,\mathcal{T}) >0$ and $\delta_I(S) = h^0(Y,\mathcal{T})-1$ by \Cref{lem:invariants}. Since $h^0(Y,\mathcal{T}) \geq  h^0(F,\O_F)-1$ by \Cref{prop:sections_wild}, it follows that $\delta_I(S)\geq h^0(F,\O_F) -2$.

(3) This follows from \cite[Proposition 2.1 (2)]{Liedtke_note} (or also \Cref{lem:invariants}).
\end{proof}

\begin{remark}
    Surfaces $S = E\times^G X$ admitting a unique wild fiber and having a singular Picard scheme $\Pic_S$ do exist. 
    
    If $\kappa(S)=0$, as explained in \cite{BMIII,BMII} (see also \cite[Theorem 1.2]{Liedtke_note}), bi-elliptic and quasi-bi-elliptic surfaces with reduced $\underline{\mathrm{Pic}}_{S}$ are precisely those with nontrivial canonical bundle. According to the classification of Mitsui \cite{Mitsui}, there exist such surfaces with exactly one wild fiber.
    
    If $\kappa(S)=1$, those are examples of Katsura and Ueno, see \cite[Example 8.1]{Katsura_Ueno}. By \cite[Proposition 2.3]{Liedtke_note}, they also provide examples with arbitrarily large $\delta_I(S)$. The new examples coming from \Cref{prop:E} are also of this form.
\end{remark}

For our isotrivial elliptic surfaces $S=E\times^G X$, the multiplicity of the fibers are entirely determined by the quotient $\pi\colon X\to Y$. More precisely, we have the following lemma:

\begin{proposition}\label{prop:stab_multiplicity}
    Assume the action of \(G\) on \(X\) is faithful.
    Let $x\in X$ and $y\in Y$ be geometric points such that $\pi(x)=y$. Then the schematic fiber $f^{-1}(y)$ has multiplicity
    \[
    m_y =\vert \mathrm{Stab}_G(x) \vert,
    \]
    where $\vert \mathrm{Stab}_G(x) \vert$ is the order of the subgroup scheme $\mathrm{Stab}_G(x)\subset G$.
\end{proposition}

\begin{proof}
    Let $H=\mathrm{Stab}_G(x)$ and $E'=E/H$. Then the morphism
    \[
    \begin{array}{ccc}
         E' & \longrightarrow & S\\
         q \mod H& \longmapsto & (q,x) \mod G
    \end{array}
    \]
    is a closed immersion which identifies $E'$ with $(f^{-1}(y))_{red}$. In particular, we have that $f^{-1}(y)=m_yE'$ as Cartier divisors, and the restriction of $h\colon S\to E/G$ to the subscheme $(f^{-1}(y))_{red}$ corresponds to the quotient of $E'$ by $G/H$. Therefore, denoting by $q_0$ the neutral point of the elliptic curve $E/G$ and $F=h^{-1}(q_0)$, we have that
    \[
    ((f^{-1}(y))_{red}\cdot F) = \vert G/H\vert,
    \]
    where the left-hand side is the intersection number of $(f^{-1}(y))_{red}$ and $F$.

    Moreover, the fiber $F$ is isomorphic to $X$, which is reduced. 
    Explicitly, this isomorphism is given by the closed immersion
    \[
    \begin{array}{ccc}
         X & \longrightarrow & S\\
         z & \longmapsto & (q_0,z) \mod G.
    \end{array}
    \]
    Via this identification, the restriction $f_{\vert F}$ coincides with the quotient morphism $\pi\colon X\to Y$. Since $\pi$ is finite morphism of degree $\vert G \vert$, this shows the equality
    \[
    (f^{-1}(y)\cdot F) = \vert G \vert,
    \]
    where the left-hand side is the intersection number of the subschemes $f^{-1}(y)$ and $F$.

    Using that $F$ is reduced and that $f^{-1}(y)=m_yE'$, we have 
    \[
    (f^{-1}(y)\cdot F) = m_y((f^{-1}(y))_{red}\cdot F).
    \]
    Thus, the above computation shows that $m_y=\vert H\vert$, and this concludes the proof.
\end{proof}

As a immediate consequence of Propositions \ref{prop:one_wild_fiber} and \ref{prop:stab_multiplicity}, we have:

\begin{corollary}\label{coro:free_action_no_multiple}
    If $G$ acts freely on the projective $G$-normal curve $X$, then the isotrivial elliptic surface 
    \[
    S=E\times^G X \to Y=X/G
    \]
    has no multiple fibers. Moreover, the Picard scheme $\Pic_S$ is reduced.
\end{corollary}

\begin{remark}
    If $\pi\colon X\to Y$ is an \emph{unramified} Artin-Schreier cover, then $\Pic_S$ is reduced. Indeed, $\pi$ is a $(\mathbb{Z}/p\mathbb{Z})$-torsor. By \Cref{coro:free_action_no_multiple}, the elliptic fibration $f\colon S\to Y$ has no multiple fiber and $\Pic_S$ is reduced.
\end{remark}

\subsection{Generalized Igusa surfaces from unipotent group actions}

If $G\subset E$ is a diagonalizable finite subgroup scheme, then $\Pic_S$ is reduced by \Cref{thm:diag_reduced_pic}. However, $E$ also contains unipotent finite subgroup schemes. For example, the $p$-torsion subgroup of the elliptic curve $E$ is
\[
E[p]
=
\begin{cases}
   (\mathbb{Z}/p\mathbb{Z} ) \times \boldsymbol{\mu}_p \text{ if } E \text{ is ordinary},\\
   I_{1,1} \text{ if } E \text{ is supersingular}.\\
\end{cases}
\]
In the supersingular case, $I_{1,1}$ is a local-local finite commutative group scheme of order $p^2$ which sits in a nontrivial extension of $\aalpha_p$ by $\aalpha_p$. For arbitrary finite subgroup schemes of $E$, we have the following lemma:

\begin{lemma}\label{group decomp}
    If \(G\) is a finite commutative group scheme over a perfect field of characteristic \(p\), then it is a product of the form
    \[G^{lin}\times G^{unip}\]
    where \(G^{lin}\) is linearly reductive and \(G^{unip}\) is unipotent.
\end{lemma}
\begin{proof}
    By \cite[Theorem 15.5]{FiniteGruopSchemesPink}, there is a unique decomposition of $G$ as
    \[G\cong G_{rr}\times G_{lr}\times G_{rl}\times G_{ll}\]
    where \(G_{xy}\subset G\) is the maximal subgroup with property \(x\) and Cartier dual with property \(y\), and where \(r\) stands for reduced and \(l\) stands for local. 
    When \(G\) and its Cartier dual are reduced, then \(G\) is of order prime to \(p\) and therefore linearly reductive, similarly, if \(G\) is local with reduced Cartier dual, then \(G\) is of multiplicative type (by definition as it is the Cartier dual of a finite abstract group).  If \(G\) is reduced but has a local Cartier dual, then its dual is of order \(p^n\), and hence the order of \(G\) is also of order \(p^n\). Thus \(G\) is an étale \(p\)-group, hence it is unipotent. 
    Finally, if \(G\) is local with local Cartier dual, then by \cite[Proposition 22.5]{FiniteGruopSchemesPink}, \(G\) is a subgroup scheme of a direct sum of groups of truncated Witt vectors. As all groups of truncted Witt vectors are unipotent and subgroups of unipotent groups are unipotent, \(G\) is unipotent as well. The desired statement then follows by setting \(G^{lin}:=G_{rr}\times G_{lr}\) and \(G^{unip}= G_{rl}\times G_{ll}\).
\end{proof}

\begin{proposition}\label{rational quotient case}
    Let \(S=E\times^GX\) and assume that \(Y=X/G\cong \mathbb{P}^1\). Then \(\Pic_S\) is singular if and only if the linearly reductive factor $G^{lin} \subset G$ acts on \(H^0(X,\omega_X)\) with nontrivial invariants. Moreover, if \(X\) has arithmetic genus \(1\), then \(\Pic_S\) is singular if and only if \(G^{lin}\) acts trivially on \(H^0(X,\omega_X)\).
\end{proposition}
\begin{proof}
    By \Cref{group decomp}, we can write \(G=G^{lin}\times G^{unip}\) as a product of a linearly reductive and a unipotent group scheme. By \Cref{lem:invariants}, 
    \(\Pic_S\) is reduced iff \(\delta_I(S) = h^0(\omega_X)^G - g(Y)=0\), but in our case \(g(Y)=0\) so \(\Pic_S\) is reduced if and only if \(h^0(\omega_X)^G=0\). 

    Since $G^{lin}$ is commutative, it is diagonalizable (see e.g. \cite[Proposition 12.54]{Milne}). Restricting the \(G\)-representation \(V=H^0(X,\omega_X)\) to \(G^{lin}\) gives a direct sum decomposition  \begin{equation}\tag{$\clubsuit$}\label{eq:decomposition_V}
      V=\bigoplus_{\chi\in\mathrm{Hom}_{\mathrm{gp}}(G^{lin}, \mathbb{G}_m)}V_{\chi}  
    \end{equation}
    indexed by the characters of \(G^{lin}\) where \(V_{\chi}\) is defined by \(g\cdot v=\chi(g)v\) for all \(v\in V_{\chi}\) and all \(g\in G^{lin}\). By the direct product decomposition and the fact that \(G\) acts linearly on \(V\), for all \(g_l\in G^{lin}\), \(g_u\in G^{unip}\) and all \(v\in V_{\chi}\) we have
    \[
    g_l\cdot (g_u\cdot v)=g_u\cdot (g_l\cdot v)=g_u\cdot (\chi(g_l)v)=\chi(g_l)(g_u\cdot v).
    \]
    In other words, \(V_{\chi}\) is stable under the action of \(G^{unip}\); thus, the above decomposition \eqref{eq:decomposition_V} is valid as a \(G\)-representation and not only as a \(G^{lin}\)-representation.
    In particular, the fixed points of \(V\) are \(V_0^{G^{unip}}\), which is always nontrivial if \(V_0\neq0\). This shows the first equivalence.
    
    For the last statement, note that the dimension of \(H^0(X,\omega_X)\) is the arithmetic genus of \(X\), so if the arithmetic genus of \(X\) is \(1\), then \(h^0(\omega_X)^G\neq 0\) if and only if the action of \(G^{lin}\) on \(H^0(X,\omega_X)\) is trivial.
\end{proof}

Over a smooth projective curve of positive genus $g(Y)>0$, the condition $h^0(\omega_X)^G>0$ no longer implies that $\delta_I(S)>0$. To deal with the higher genus cases, we show that $G$ having fixed points guarantees the existence of wild fibers when \(G\) is unipotent:

\begin{proposition}\label{coro:unipotent_stab_wild}
    Assume that $G$ acts faithfully on $X$. If \(x\in X\) has a nontrivial unipotent stabilizer, then \(f^{-1}(\pi(x))\) is a wild fiber. 
\end{proposition}

\begin{proof}
    Let \(H=\mathrm{Stab}_G(x)\), \(y:=\pi(x)\) and \(F=f^{-1}(y)\).
    We need to show that \(h^0(F,\mathcal{O}_F)\geq 2\). First, note that
    \[
    H^0(F,\mathcal{O}_F)=H^0(E\times\pi^{-1}(y),\mathcal{O}_{E\times\pi^{-1}(y)})^G=H^0(\pi^{-1}(y),\mathcal{O}_{\pi^{-1}(y)})^G
    \]
    where the first equality is by the definition of \(f:S\to Y\) and descent, and the second follows from the Künneth formula.
    Note that \(\pi^{-1}(y)\) is finite over \(\Spec k\), in particular it is affine, and \(\pi^{-1}(y)=\Spec A\) for the finite \(k\)-algebra \(A=H^{0}(\pi^{-1}(y),\mathcal{O}_{\pi^{-1}(y)})\). 
    
    The action \(G\times X\to X\) restricted to \(G\times\{x\}\) induces a \(G\)-invariant closed immersion \(G/H\subset \pi^{-1}(y)=\Spec A\) defined by an ideal \(I\), which is \(G\)-invariant. Note that \(\dim_k k[G/H]=|G/H|<|G|=\dim_kA\) as the quotient \(X\to Y\) is finite flat of degree \(|G|\) by \Cref{cor:fin-flat}
    and \(H\) is assumed non-trivial; and thus, \(I\) is a non-trivial ideal. In particular, its elements are non-constant functions on \(\Spec A\) or the zero function. 

    By \Cref{group decomp}, \(G=U\times L\) where \(U\) is unipotent and \(L\) is diagonalizable, in particular \(H^0(F,\mathcal{O}_F)=A^G=(A^L)^U\).
    Since we assumed that \(H\) is unipotent, it follows that \(H\subset U\) and \(H\cap L=\{e\}\). In other words, the fiber of the finite group scheme \(\mathrm{Stab}_L\subset L\times X\to X\) over \(x\) is trivial, and all closed points in \(\pi^{-1}(y)\) have the same \(L\)-stabilizer by commutativity. Hence by Nakayama's lemma, it is trivial on all of \(\Spec A\), i.e the action of \(L\) on \(\Spec A\) is free. In particular, \(\Spec A\to \Spec A^{L}\) is an \(L\)-torsor.

    Since \(L\) is linearly reductive, \((A/I)^L=A^L/I^L\), so the ideal defining the closed immersion \(G/(H\times L)\subset \Spec A^L\), coming from descent of invariant closed subschemes along torsors, is \(I^L\subset A^L\).
    As invariant closed subschemes of \(\Spec A\) and closed subschemes of \(\Spec A^L\) are in one-to-one correspondence by descent, it follows that \(I^L\subset A^L\) is a non-trivial ideal since \(G/H\) is neither empty nor equal to \(\Spec A\). In particular, \(I^L\) is a non-trivial ideal of \(A^{L}\).
    
   Now, \((I^{L})^{U}= I^{G}\subset A^{G}=H^{0}(\pi^{-1}(y),\mathcal{O}_{\pi^{-1}(y)})^{G}\) is a non-trivial vector subspace, as \(U\) is unipotent and \(I^{L}\) is a non-zero \(U\)-representation. Furthermore, \(I^G\cap k=\{0\}\) where \(k\) is viewed as the ring of constant functions on \(\Spec A\), so \(\dim_{k}A^G\geq \dim_k(k\oplus I^G)\geq 2\).
    \end{proof}

A particularly useful special case is:
\begin{corollary}
\label{lemma unipotent wild}
    Let \(x\in X\) be a fixed point with \(y=\pi(x)\in Y\) and assume that \(G\neq\{e\}\) is unipotent and acts faithfully on $X$, then \(F=f^{-1}(y)\) is a wild fiber.
\end{corollary}
\begin{proof}
    This is a direct consequence of \Cref{coro:unipotent_stab_wild}.
\end{proof}

\begin{proof}[Proof of \Cref{thm:A}]
The statements follow from \Cref{lem:invariants}, \Cref{prop:one_wild_fiber,prop:stab_multiplicity}, and \Cref{coro:unipotent_stab_wild}.
\end{proof}

\begin{proof}[Proof of \Cref{coro:B}]
If $x_1,x_2$ are two points with distinct orbits and nontrivial stabilizers, then the fibers at those two points are wild, so $S$ is a generalized Igusa surface by \Cref{prop:one_wild_fiber}. The second statement is a special case of \Cref{rational quotient case}.
\end{proof}

\begin{remark}
    In Igusa's construction, we have that $p=2$ and $X$ is an elliptic curve. Then the quotient $\pi\colon X\to \P^1$ by the involution $x\mapsto -x$ is an Artin-Schreier cover. Notice that it is ramified at the two-torsion points, so the resulting elliptic fibration $f\colon S=E\times^{\mathbb{Z}/2\mathbb{Z}} X\to \P^1$ has one or two multiple fibers by \Cref{prop:stab_multiplicity}, depending on whether $X$ is supersingular or ordinary. One can see that $\Pic_S$ is singular with \Cref{coro:B}.
\end{remark}

\section{Constructing generalized Igusa surfaces}\label{section:applications}

\subsection{Generalized Igusa surfaces via Artin-Schreier curves}\label{subsection:generalized_Igusa_Z/p}

~\medskip

Let us recall some definitions and properties from \cite[Chapters IV and X]{Serre_local_fields} and \cite[\S 7]{GarnekI}. A field extension $L/K$ is called an \emph{Artin-Schreier extension} if $L=K(\beta)$ for some $\beta \in L\setminus K$ such that $\beta^p-\beta \in K$. When \(K=k(Y)\) is the function field of a curve \(Y\), every (potentially ramified) \((\mathbb{Z}/p\mathbb{Z})\)-Galois cover \(\pi: X\to Y\) of \(Y\) can be described by an Artin-Schreier extension of \(K\). Given a point \(y\in Y\) of such a Galois cover, the equation 
\[
x^p-x=\alpha
\]
is of \emph{local standard form} at \(y\) if \(\alpha\) is regular at \(y\) or has a pole order coprime to \(p\), the latter happens only if \(y\) is a branch point for the cover. When \(y\) is a branch point, the pole order of a function \(\alpha\), giving a local standard form description of \(\pi\colon X\to Y\), is independent of the choice of such representation. More precisely, the $i$-th ramification group of $\pi$ at $y=\pi(x)$ is 
\[
G_{y,i} = \{\sigma \in \mathbb{Z}/p\mathbb{Z}~\vert ~ \sigma \cdot f \equiv f ~\mathrm{mod}(\mathfrak{m}_x^{i+1}) ~\text{for every }f\in \O_{X,x}\},
\]
where $\mathfrak{m}_x$ denotes the unique maximal ideal of $\O_{X,x}$. 
In particular, $G_{y,i}$ is trivial or equals  $\mathbb{Z}/p\mathbb{Z}$, and the \emph{ramification jump} at $y$ is the integer 
\[
m_y=\max\{i~\vert ~G_{y,i}\neq 0\}.
\]
The integer $m_y$ coincides with the pole order of $\alpha$ at \(y\) (see e.g. \cite[\S 7.1]{GarnekI}).

\medskip

\begin{remark}
The Igusa surface has already been generalized by Katsura and Ueno in \cite[Example 8.1]{Katsura_Ueno}: they consider the contracted product $S=E\times^{\mathbb{Z}/p\mathbb{Z}} X$, where $E$ is an ordinary elliptic curve and $X$ is an Artin-Schreier curve over $\P^1$ ramified at \(\infty\). For these surfaces, Liedtke shows that taking the pole at infinity with sufficiently large multiplicity $m$, the Picard scheme is singular (see \cite[Proposition 2.3]{Liedtke_note}).     
\end{remark}

Let $E$ be an ordinary elliptic curve and let $\pi\colon X\to Y$ be an Artin-Schreier cover. The $p$-torsion subgroup $E[p]\subset E$ contains a copy of $\mathbb{Z}/p\mathbb{Z}$, acting on $E$ by translations, and on $X$ via the Artin-Schreier action. In this section, we consider the contracted product of $E$ and $X$ with respect to the induced diagonal action
\[
    S = E\times^{\mathbb{Z}/p\mathbb{Z}} X.
\]
If $\pi$ is ramified at two points, then the elliptic surface $S=E\times^{\mathbb{Z}/p\mathbb{Z}} X$ has a singular Picard scheme by \Cref{coro:B}. A more detailed analysis is required to deal with the case where $\pi$ is ramified at a single point. The following example, which is due to Katsura and Ueno (see \cite[Example 4.7]{Katsura_Ueno}), shows that the equality $\delta_I(S)=0$ can occur. In fact, we will see in \Cref{thm:C} that it is essentially the only example of a contracted product $E\times^{\mathbb{Z}/p\mathbb{Z}} X$ with a ramified Artin-Schreier cover $X\to Y$ and which has a reduced Picard scheme.

\begin{example}
    Let $E$ be an ordinary elliptic curve, so that $\mathbb{Z}/p\mathbb{Z}\subset E[p]$. Consider $\mathbb{Z}/p\mathbb{Z}$ acting by $x\mapsto x+1$ on $\mathbb{A}^1\subset \P^1$ and fixing the point at infinity. This yields an Artin-Schreier cover $X=\P^1\to \P^1=Y$ with an affine equation
    \[
    x^p-x = y.
    \]
    In particular, the point at infinity is a single pole of order $m=1$. Then the induced diagonal action yields a contracted product $S=E\times^{\mathbb{Z}/p\mathbb{Z}} \P^1$, that is a geometrically ruled surface. Thus, $\Pic_S$ is reduced. This can also be seen from the fact that \(H^0(\mathbb{P}^1,\omega_{\mathbb{P}^1})=0\), so \(\delta_I(S)=h^0(\mathbb{P}^1,\omega_{\mathbb{P}^1})^G-h^0(\mathbb{P}^1,\omega_{\mathbb{P}^1})=0-0=0\). 
\end{example}

By \Cref{lem:invariants}, it suffices to study the invariants of the action of $\mathbb{Z}/p\mathbb{Z}$ on $H^0(X,\omega_X)$. To this end, we use a result of Valentini and Madan, see \cite[Theorem 1]{Valentini_Madan}. Their theorem is stated in a much more general form, as they study towers of Artin-Schreier extensions. For a single Artin-Schreier cover $X\to Y$ we have:

\begin{remark}
    Let $\pi\colon X\to Y$ be an Artin-Schreier cover ramified at the points $\{p_i\}_{i=1,\ldots,s}$ in $X$. Denoting $q_i=\pi(p_i)$, we have
    \[
    m_{i} = \Phi(i,1),
    \]
    where $m_{i}$ is the ramification jump at $q_i$ and $\Phi(i,1)$ is the positive integer defined in \cite[p. 110]{Valentini_Madan}. 

    Indeed, using their notations, we have in our situation: $j=1$ since we work with a single Artin-Schreier cover instead of arbitrary towers of extensions, $m=1$ since every branched point $q_i\in Y$ is the image of unique ramification point $p_i\in X$, and $b_1=\alpha\in k(Y)$ is a rational function such that the equation $x^p-x=\alpha$ is of local standard form at every $q_i$ (see the properties (1) and (2) of \cite[p. 108]{Valentini_Madan}). Therefore,
    \[
    \Phi(i,1) = -\nu(i,0,1,\alpha)
    \]
    is the pole order of $\alpha$ at $q_i$, as $\nu(i,0,1,\alpha)$ denotes the normalized valuation at $q_i$. This coincides with the ramification jump $m_{i}$ as noticed in \cite[\S 7.1]{GarnekI}.
\end{remark}

Analogously as in \cite[Example 28]{Karanikolopoulos_Kontogeorgis}, we apply their theorem to our extension $k(Y) \subset k(X)$ and we get:

\begin{theorem}[{\cite[Theorem 1]{Valentini_Madan}}]\label{thm:Valentini_Madan}
    Let $G=\mathbb{Z}/p\mathbb{Z}$ and let $\pi\colon X\to Y=X/G$ be a ramified Artin-Schreier cover. Let $\{q_i\}_{i=1,\ldots,s}$ be the branch points of the cover, each of ramification jump $m_i>0$. Then the following hold:
    \begin{enumerate}
        \item The regular representation of $G$ occurs $g(Y)$ times in the representation of $G$ on $H^0(X,\omega_X)$.
        \item For $d=1,\ldots,p-1$, the indecomposable representation of degree $d$ occurs $\Gamma_{d-1} - \Gamma_d +a_d$ times in the representation of $G$ on $H^0(X,\omega_X)$, where 
        \[
        a_d=
        \begin{cases}
            -1 & \text{ if }d=p-1 \\
            0  &\text{ if }d<p-1,
        \end{cases}
        \]
        and for $d=0,\ldots,p-1$,
        \[
        \Gamma_d = \sum_{i=1}^s \nu_{id},
       ~\text{and }~
        \nu_{id}=\left\lfloor \frac{(m_i+1)(p-1)-dm_i}{p} \right\rfloor .
        \]
    \end{enumerate}
    
\end{theorem}

Applying the above theorem to the smooth projective surfaces $S=E\times^{\mathbb{Z}/p\mathbb{Z}} X$, we obtain:

\begin{proposition}\label{prop:Delta_bound_AS}
    Let $\pi\colon X\to Y$ be an Artin-Schreier cover as in  \Cref{thm:Valentini_Madan}. Let $E$ be an ordinary elliptic curve and let $S=E\times^G X$ be the induced contracted product. Then 
    \[
    \delta_I(S) =  \Gamma_0 -1 = \left( \sum_{i=1}^s \left\lfloor \frac{(m_i+1)(p-1)}{p} \right\rfloor \right) -1,
    \]
    and $\delta_I(S)=0$ if and only if $\Gamma_0=1$. 
\end{proposition}

\begin{proof}
     For each regular representation $k[G]\subset H^0(X,\omega_X)$, we have $k[G]^G=k$. By \Cref{thm:Valentini_Madan}, the regular representation occurs $g(Y)$ times. Since
        \[ \dim\left(\bigoplus_{i=1}^{g(Y)}k[G] \right)^G = g(Y),
        \]
        it follows by \Cref{lem:invariants} that $\delta_I(S)=h^0(X,\omega_X)^G - g(Y)>0$ if and only if there is at least one indecomposable representation of degree $d=1,\ldots,p-1$ that occurs in the representation of $G$ on $H^0(X,\omega_X)$. Else, $\delta_I(S)=0$.

        For each indecomposable representation $V$ of degree $d=1,\ldots,p-1$, we have as well that $\dim(V^G)= 1$ (see \cite[Theorem 6.1.2, Example 6.1.3]{P.Webb}). Therefore, again by \Cref{thm:Valentini_Madan},
        \[
        \delta_I(S)=h^0(\omega_X)^G - g(Y) = \sum_{d=1}^{p-1} \big(\Gamma_{d-1} - \Gamma_d +a_d\big).
        \]
        This is a telescoping sum, which equals the non-negative integer $\Gamma_{0} - \Gamma_{p-1} -1$. Moreover, notice that for every $p\geq 2$, we have $\Gamma_{p-1}=0$ and 
        \[
        \Gamma_0 = \sum_{i=1}^s \left\lfloor \frac{(m_i+1)(p-1)}{p} \right\rfloor .
        \]
        This concludes the proof, as $\Gamma_0=1$ is equivalent to the fact that none of the indecomposable representations of degree $d=1,\ldots,p-1$ occurs. 
\end{proof}

\begin{corollary}\label{coro:Delta_bound_AS}
    Let $E$ be an ordinary elliptic curve and let $\pi\colon X\to Y$ be an Artin-Schreier cover as in \Cref{thm:Valentini_Madan}, and let $S=E\times^G X$ be the induced contracted product. We denote by
    \[
    m=\sum_{i=1}^s m_i,
    \]
    which is a positive integer. Then 
    $\delta_I(S)=0$ if and only if $s=m=1$.
\end{corollary}

\begin{proof}
    For each $i=1,\ldots,s$, we have that
    \[
    \left\lfloor \frac{(m_i+1)(p-1)}{p} \right\rfloor \geq 1,
    \]
    as $m_i >0$. Hence, $\Gamma_0=1$ if and only if $s=1$ and $m=1$. Otherwise, $\Gamma_0\geq 2$, so $\delta_I(S)>0$ as well by \Cref{prop:Delta_bound_AS}.
\end{proof}

We are now ready to prove:

\begin{proof}[Proof of \Cref{thm:C}]
    Let \(\pi:X\to Y\) be an Artin-Schreier cover, \(E\) an ordinary elliptic curve and \(S=E\times^{(\mathbb{Z}/p\mathbb{Z})} X\) the associated contracted product. If $\pi$ is unramified, we are done by \Cref{coro:free_action_no_multiple}. Indeed, ramification points of \(\pi\) correspond to \(\mathbb{Z}/p\mathbb{Z}\)-fixed points of \(X\). Else, the computation of \(\delta_I(S)\) is \Cref{prop:Delta_bound_AS} and the characterization of when \(S\) is generalized Igusa follows from \Cref{coro:Delta_bound_AS}.
\end{proof}    

These generalized Igusa surfaces provide new examples where \(\delta_I(S)\) is arbitrarily large. The existence of such surfaces was pointed out in \cite[Proposition 2.3]{Liedtke_note}, using the construction in \cite[Example 8.1]{Katsura_Ueno} of certain \(\mathbb{Z}/p^n\mathbb{Z}\)-covers of \(\mathbb{P}^1\). Moreover, these examples show that \(\delta_I(S)\) can be arbitrarily large for elliptic fibrations \(f:S\to C\) over any smooth curve, not just over \(\mathbb{P}^1\).

\subsection{Generalized Igusa surfaces via Hilario-Schröer's $\aalpha_p$-normal curves}\label{subsection:generalized_Igusa_alpha_p}

\begin{remark}
    Let us first recall the following classical construction in $p=\mathrm{char}(k)\in \{2,3\}$, that can be found in \cite[\S 2]{BMIII}. Let $E$ be a supersingular elliptic curve (in particular, $\aalpha_p\subset E[p]$) and $C_0$ be the projective rational cuspidal curve, whose affine equation is
    \[  
    x^3 = y^2.
    \]
    A \emph{quasi-bielliptic surface} is a smooth projective surface of the form $E\times^G C_0$, where $G\subset E$ is a finite group scheme, and $C_0$ is $G$-normal (see \cite[Corollary 4.8]{G-normal-curves}). A list classifying such finite group schemes $G\subset E$ is given \cite[\S 2, p. 214]{BMIII}: in particular, $C_0$ is $\aalpha_p$-normal. Moreover, by \cite[Theorem 1.2]{Liedtke_note}, the induced smooth projective surface $E\times^G C_0$ may have a singular Picard scheme.
\end{remark}

Instead of considering quotients by a constant group, we are now interested in the case where $G$ is \emph{infinitesimal}.
In this section, we aim to exploit some constructions due to Hilario and Schröer in \cite{Generalized-quasielliptic}, which hold for $p=\mathrm{char}(k)>0$. More precisely, they study a family $X_{p,n}$ of projective $U_n$-normal curves, which are generalizations of the above rational cuspidal curve $C_0$, and where $U_n$ is an infinitesimal group scheme of order $p^{n(n+1)/2}$ defined in \cite[Section \S 4]{Generalized-quasielliptic}. We focus on the case where $n=1$, so that $U_1=\aalpha_p$. Among their various results on the geometry of these curves, let us summarize the information that we need in the following.

\begin{theorem}[{\cite[Corollary 6.2, Theorem 9.4]{Generalized-quasielliptic}}]\label{Hilario Schröer}
    Let $X_p$ be the projective curve defined by 
    \[\Spec k[T^p,T^{p-1}]\cup \Spec k[T^{-1}]\]
    glued along \(\{T^{-1},T\neq0\}\).  Then the following hold:
    \begin{enumerate}
        \item There is an action of $\aalpha_p$ such that \(X_p\) is \(\boldsymbol{\alpha}_p\)-normal.
        \item \(h^1(X_p,\O_{X_p})= \frac{1}{2}\left(p^2-3p+2\right)\).
    \end{enumerate} 
\end{theorem}

Notice that $X_2 \cong \P^1$, and $X_3 \cong C_0$ is the above projective rational cuspidal curve. For $p>0$, $X_p$ comes with an action of \(\boldsymbol{\alpha}_p\), as a coaction on \(k[T^{-1}]\) defined by
\begin{equation}\label{eq:coaction_Xp} \tag{$\dagger$}
    T^{-1}\mapsto T^{-1}+sT^{-p}
\end{equation}
Indeed, this is the special case $n=1$ of the action defined in \cite[\S 5]{Generalized-quasielliptic}, the corresponding notation being $\lambda_0=1$, and $\lambda_1=s$ and $\mu = T^{-1}$.

Moreover, the normalization of \(X_p\) is \(\mathbb{P}^1\), and the composition
\[\mathbb{P}^1\to X_p\to X_p/{\boldsymbol{\alpha}_p}\]
is surjective. By \cite[Corollary 4.16]{G-normal-curves}, \(X_p/\boldsymbol{\alpha}_p\) is smooth. Hence, it follows by Lüroth's theorem
that \(X_p/\boldsymbol{\alpha}_p\cong \mathbb{P}^1\). As for the case $G=\mathbb{Z}/p\mathbb{Z}$, we consider the contracted product
\[
S_p = E\times^{\aalpha_p} X_p,
\]
where $E$ is now a supersingular elliptic curve.
We can already notice that for $p\geq 3$, we have that $p_a(X_p)\geq 1$; thus, it follows from \Cref{rational quotient case} that $S_p$ is a generalized Igusa surface. We make this statement more precise in \Cref{prop:E}, that we prove here:

\begin{proof}[Proof of \Cref{prop:E}]
    Notice that $S_2$ is a geometrically ruled surface; hence, $\Pic_{S_2}$ is reduced by \cite[Theorem 1.1]{Liedtke_note}. Notice also the surface $S_3$ is one of the quasi-bielliptic surfaces considered by Bombieri-Mumford in \cite[\S 2]{BMIII}. When $p\geq 5$, the fact that $S_p$ has Kodaira dimension one follows readily from $p_a(X_p)\geq 2$ (\Cref{Hilario Schröer}) together with \cite[Theorem A]{fong_maccan}. 

    By \Cref{lem:invariants}, since $X_p/\aalpha_p \cong \P^1$, it remains to show that $h^0(\omega_{X_p})^{\aalpha_p}=p-2$. As $X_p\cong \P^1$ if $p=2$, we can assume from now on that $p\geq 3$.
    Let us consider the projective model of $X_p$ described in \cite[\S 7, (7.2)]{Generalized-quasielliptic}: there is a natural isomorphism
    \[
    X_p \longrightarrow  C = \{U^p-V^{p-1}Z=0\} \subset \mathbb P^2_{U,V,Z}, \qquad T^{p-1} \mapsto U/Z, \quad T^p \mapsto V/Z, \quad T^{-1} \mapsto U/V.
    \]
    Using this projective model, we have $\omega_{X_p} = \O_{X_p}(p-3)$ by \cite[Proposition 7.1]{Generalized-quasielliptic}. Taking the global sections in the following short exact sequence
    \[
    0 \to \O_{\P^2}(-3) \to \O_{\P^2}(p-3) \to \O_{X_p}(p-3) \to 0,
    \]
    we deduce that $h^0(\omega_{X_p}) \cong h^0(\P^2,\O_{\P^2}(p-3))$, which is isomorphic to $\mathrm{Sym}^{p-3}(\mathfrak{g})$, where $\mathfrak{g}=H^0(X,\Theta_X)$ is the $3$-dimensional vector space of global sections of the tangent sheaf $\Theta_X$, with basis $\{Z,V,U\}$ (see \cite[Proposition 7.2]{Generalized-quasielliptic}).

    Via the above identification, it follows from the proof of \cite[Theorem 8.1, case $n=1$]{Generalized-quasielliptic} (see page 21 of loc.cit.) that the infinitesimal group scheme $\aalpha_p$ acts on $\mathfrak{g}$ by the matrices
    \[
    \begin{pmatrix}
        1 & 0 & c \\ 0 & 1 & 0 \\ 0 & 0 & 1
    \end{pmatrix}, ~\text{where }c^p=0
    \]
    Equivalently, the coaction is described by
    \[
        Z \mapsto Z,\quad 
        V \mapsto V, \quad
        U \mapsto U+cZ.
    \]
    Writing an homogeneous polynomial $P$ of degree $p-3$ as
    \[
    P = \sum_{i=0}^{p-3} a_{i}(Z,V)U^{i}
    \]
    where $a_{i}(Z,V)$ has degree $p-3-i$, we see that $P$ is $\aalpha_p$-invariant if and only if $P=a_0(Z,V)$ is independent of $U$. This implies that $(\mathrm{Sym}^{p-3}(\mathfrak{g}))^{\alpha_p}$ is the linear subspace generated by $Z$ and $V$, which has dimension $p-2$.
\end{proof}

\begin{remark}
    The surfaces \(S_p\), when \(p\geq 5\), can be distinguished from the previously known examples of generalized Igusa surfaces via their automorphism groups. Indeed, to our knowledge, all the previous Kodaira dimension one examples have \((\Aut_{S}^{0})_{\mathrm{red}}\) as an ordinary elliptic curve while \((\Aut_{S_p}^{0})_{\mathrm{red}}\) are supersingular elliptic curves for all \(p\) (see \Cref{prop:Rosenlicht_application}).
\end{remark}

\section{Singular Picard schemes of (quasi)-bielliptic surfaces}
\label{sec: picinf}
Let us recall that $\Pic^0_Z$ is a group scheme of finite type over $k$, for any projective variety $Z$. The base field being algebraically closed implies that the reduced subscheme of any algebraic group $H$ forms a closed subgroup (in general this is true over perfect base fields), which is also normal as \(\Pic_Z^0\) is commutative. We denote its quotient, which is an infinitesimal group scheme, as
\[
H_{\mathrm{inf}} := H/H_{\mathrm{red}}
\]
Let $S$ be a generalized Igusa surface. We aim to compute the group scheme structure of the infinitesimal quotient $(\Pic_S^0)_{\mathrm{inf}}$, for all such $S$ that have Kodaira dimension zero.  Among the latter, generalized Igusa surfaces do occur, only for $p=2,3$ and only among (quasi)-bielliptic ones. We first recall some general properties of Picard schemes of quotients, that allow us to reduce to computing fixed points of the Picard scheme of a \(G\)-normal curve, and we describe the Picard scheme of the cuspidal curve. Then, we describe the group scheme structure of $(\Pic_S^0)_{\mathrm{inf}}$, for $S$ quasi-bielliptic, and then we do the same for $S$ bielliptic.

\subsection{Picard schemes of quotients}
Let us start with recalling the following, which gives us a tool to compute Picard schemes of quotients; see \cite[Theorem 2.1]{Jensen}.

\begin{theorem}[{\cite[Theorem 2.1]{Jensen}}]\label{thm:Jensen} 
    Let $G$ be a finite commutative group scheme acting freely on a projective variety $W$ and assume that the quotient $S := W/G$ exists. Assume moreover that: either $G$ is constant and cyclic, or $G$ is constant and $D(G)$ is infinitesimal or vice-versa, or finally that $G$ is infinitesimal of height one. Then there is a short exact sequence
    \[
    0 \longrightarrow D(G) \longrightarrow \underline{\mathrm{Pic}}_{S} \longrightarrow (\underline{\mathrm{Pic}}_W)^G \longrightarrow 0
    \]
    of commutative group schemes, induced by the pull-back map.
\end{theorem}

Following the same arguments as in \cite[Example 2.2]{Jensen}, we restrict the short exact sequence to the torsion component $\Pic_S^{\tau}\subset \Pic_S$, which consists in line bundles whose classes are torsion in $\NS_S=\Pic_S/\Pic^0_S$, to obtain the following.

\begin{corollary}\label{cor:Jensen}
Let $S$ be a generalized Igusa surface. Assume that $G$ is as in \Cref{thm:Jensen}, and moreover that there is a fixed point for the $G$-action on \(X\). Then
    \[
    \Pic_S^{\tau}  = \Pic_{E/G}^{0} \times \, (\Pic^0_{X})^G.
    \]
   In particular, there is an isomorphism 
   $(\Pic_S^0)_{\mathrm{inf}} \simeq ((\Pic_X^0)^G)_{\mathrm{inf}}$.
\end{corollary}

\begin{proof}
    Restricting the short exact sequence in \Cref{thm:Jensen} to the torsion components gives
    \[
    0 \to D(G) \to \Pic_S^\tau \to (\Pic_{E\times X}^\tau)^G \to 0. 
    \]
    Let $p_E\colon E\times X \to E$ and $h\colon S\to E/G$ be the natural projections. Then we have the following commutative diagram
    \[
    \begin{tikzcd}
        & & 0 & 0 & \\
        & & Q \ar[r,"\cong"] \ar[u] & (\Pic_X^\tau)^G \ar[u] & \\
        0 \ar[r] & D(G) \ar[r] & \Pic_S^\tau \ar[r] \ar[u] &  (\Pic_{E\times X}^\tau)^G \ar[r] \ar[u] & 0 \\
        0 \ar[r] & D(G) \ar[r]\ar[u,equal] & \Pic_{E/G}^\tau \ar[r] \ar[u,"h^*"] &  (\Pic_{E}^\tau)^G \ar[r] \ar[u,"p_E^*"] & 0 \\
        & & 0 \ar[u] & 0 \ar[u] &,
    \end{tikzcd}
    \]
    where $Q$ denotes the quotient of $\Pic_S^\tau$ by $h^*(\Pic_{E/G}^\tau)$, and the two vertical and the two horizontal lines are short exact sequences.

    We claim that the injection \(h^*\colon \Pic^{\tau}_{E/G}\to \Pic^{\tau}_{S}\) admits a retraction, and hence that \(\Pic^{\tau}_S\cong\Pic^{\tau}_{E/G}\times Q\), since the category of commutative group schemes is abelian and the existence of a retraction of a short exact sequence implies the existence of a section. Since there is a fixed point for the action of \(G\) on \(X\) by assumption, there is a point \(p\in Y=X/G\) such that \(f^{-1}(p)_{\mathrm{red}}\cong E/G\), and the closed immersion \(f^{-1}(p)_{\mathrm{red}}\cong E/G\to S\) composed with \(h\colon S \to E/G\) is easily seen to be the identity, that is \(h\colon S\to E/G\) admits a section. Pullback along this section yields the desired retraction, and hence also the splitting. Note that the splitting is canonical precisely when \(X\) has a unique fixed point. To conclude, since $E/G$ and $X$ are integral projective curves, the torsion components of their Picard schemes coincide with their connected component of identity.

    To prove the second assertion, it suffices to consider the quotient by the reduced subscheme on both sides, together with the fact that $\Pic^0_{E/G}$ is the dual elliptic curve of $E/G$.
\end{proof}

\begin{remark}
\label{arithmetic genus one}
    Assume that $\kappa(S)=0$. By \cite{fong_maccan}, $X$ is a curve with arithmetic genus one and $S$ is either an abelian surface, or a (quasi)-bielliptic surface. If $S$ is abelian, then $\Pic_S$ is reduced, so we may exclude this case. Bielliptic surfaces arise when $X$ is an elliptic curve on which $G$ does not act only by translations, see \cite[Section 3]{BMII}. Quasi-bielliptic surfaces arise when $X$ is a rational curve with a single ordinary cusp (see e.g. \cite[\S 2]{BMIII}). In both cases, $\Pic_S$ may be singular only in characteristics 2 and 3, see \cite[Theorem 1.2]{Liedtke_note}. 

    Notice that for $X$ as above, $\omega_X$ is trivial and $Y=\P^1$. Hence, $H^0(\omega_X)=k$ and $\Pic_S$ is singular if and only if $h^0(\omega_X)^G=1$ by \Cref{lem:invariants}. The latter condition is automatically fulfilled if $G$ is unipotent, as any one-dimensional representation of $G$ is trivial.
\end{remark}

To apply this in the quasi-bielliptic examples, we compute \(\Pic_{X_n}^0\) where 
\[
X_n = \Spec k[T^n,T^{n-1}] \cup \Spec k[T^{-1}]
\]
is the curve as defined in \Cref{Hilario Schröer}, but, without requiring that \(n=p\) is the characteristic of the ground field $k$.
While we will actually only carry out the full computations for surfaces defined via \(X_3\) in characteristic \(2\) and \(3\), we include the general computation of the Picard schemes for future reference.

Most of the following proposition can alternatively be proven rather straightforwardly using the conductor square in \cite[Section 5, formula (5.2)]{Generalized-quasielliptic}. However, their proof does not lend itself well to computations of the induced action of $\aalpha_p$ on \(\Pic_{X_p}^0\), since it crucially uses the normalization, and the action on \(X_p\) does not lift. 

\begin{proposition}\label{lem: picardXp}
    Let \(X_n=\Spec k[T^n,T^{n-1}]\cup\Spec k[T^{-1}]\) glued as above with \(n\geq 3\). 
    For any \(k\)-algebra \(R\), 
    \[\Pic_{X_n}(R)\cong R[T^{-1}]^{\times}\setminus R((T))^{\times}/R[[T^n,T^{n-1}]]^{\times}\]
    as abelian groups with \(\Pic^0_{X_n}(R)\) being the subset of classes represented by Laurent series of degree-valuation \(0\).
    In particular, the following hold:
    \begin{enumerate}
        \item \(\Pic_{X_3}^0\cong \mathbb{G}_a\) for any characteristic.
        \item If \(\mathrm{char} (k)=p=n\), then \(\Pic_{X_p}^0\cong \mathbb{G}_a^{\frac{(p-1)(p-2)}{2}}\). Moreover, in this case, the assignment \[T^{-1}\mapsto T^{-1}+sT^{-p}\] defines an action of \(\aalpha_p\) on \(R((T))\) by continuous \(R\)-algebra automorphisms, for any \(s\in \aalpha_p(R)\). The subgroups \(R[T^{-1}]^{\times}\) and \(R[[T^p,T^{p-1}]]^{\times}\) are invariant; this defines an action on the double quotient, recovering the natural action of \(\aalpha_p\) on \(\Pic_{X_p}\) via pullback. 
    \end{enumerate}
\end{proposition}
\begin{proof}
    Let \(s \colon \Spec k\to X_n\) be the unique singular section and let \(\nu\colon \mathbb{P}^1\to X_n\) be the normalization. Note that \(s\colon \Spec k\to X_n\) naturally factors through \(\nu\colon \mathbb{P}^1\to X_n\), we denote the lift \(\sigma \colon \Spec k\to \mathbb{P}^1\).
    By \cite[Tag 0D28]{stacks-project}, we have a commutative diagram 
    \[\begin{tikzcd}
	0 & {\mathrm{Pic}(R)} & {\mathrm{Pic}(X_n\times \Spec R)} & {\Pic_{X_n}(R)} & 0 \\
	0 & {\mathrm{Pic}(R)} & {\mathrm{Pic}(\mathbb{P}^1\times \Spec R)} & {\Pic_{\mathbb{P}^1}(R)} & 0
	\arrow[from=1-1, to=1-2]
	\arrow[from=1-2, to=1-3]
	\arrow["{=}", from=1-2, to=2-2]
	\arrow[from=1-3, to=1-4]
	\arrow["{\nu^*}", from=1-3, to=2-3]
	\arrow[from=1-4, to=1-5]
	\arrow["{\nu^*}", from=1-4, to=2-4]
	\arrow[from=2-1, to=2-2]
	\arrow[from=2-2, to=2-3]
	\arrow[from=2-3, to=2-4]
	\arrow[from=2-4, to=2-5]
\end{tikzcd}\]
    of abelian groups, where the rows are short exact sequences, split by the pullbacks along the induced sections \(s_R^*\colon \mathrm{Pic}(X_n\times\Spec R)\to \mathrm{Pic}(R)\) and \(\sigma_R^*\colon \mathrm{Pic}(\mathbb{P}^1\times\Spec R)\to \mathrm{Pic}(R)\).
    We claim that the kernel of \(s_R^*\colon \mathrm{Pic}(X_n\times\Spec R)\to \mathrm{Pic}(R)\), that is \(\Pic_{X_n}(R)\) is equal to the subgroup \(K(R)\subset \mathrm{Pic}(X_n\times \Spec R)\), defined by 
    \[K(R)=\{{L\in \mathrm{Pic}(X_n\times\Spec R)| L_{\vert  U_R}\cong\mathcal{O}_{U_R}\text{ and }L_{\vert V_{R}}\cong \mathcal{O}_{V_R}}\}\]
    where 
    \[U_R=\Spec R[T^{-1}]\cong \mathbb{A}^1_T\]
    and 
    \[V_R=\Spec R[[T^n,T^{n-1}]].\]
    First, note that \(K(R)\subset\ker(s_R^*)\), indeed \(s_R:\Spec R\to X_n\times \Spec R\) factors through \(V_R\) and \(s_R^*\mathcal{O}_{V_R}=\mathcal{O}_{\Spec R}\). To prove the converse inclusion, assume that \(s_R^*L\cong\mathcal{O}_{\Spec R}\). Note that \(s_R^*=\sigma_R^*\circ\nu_R^*\) and we may write \(\nu_R^*L=(M,\underline{n})\subset \mathrm{Pic}(\mathbb{P}^1\times\Spec R)\cong \mathrm{Pic}(R)\oplus\underline{\mathbb{Z}}\) using \(\Pic_{\mathbb{P}^1}\cong \underline{\mathbb{Z}}\) and the split exact sequence above. It follows by assumption that \(M\cong \mathcal{O}_R\) and the fact that \(\nu_R \colon \nu_R^{-1}(U_R)\to U_R\) is an isomorphism implies that \(L_{\vert U_R}\cong \mathcal{O}_{U_R}\). To see that \(L_{\vert V_R}\cong \O_{V_R}\), let \(\phi_0\colon R\to s^*L\) be an isomorphism of \(R\)-modules, which exists by assumption. Since the Jacobson radical of \(R[[T^{n-1},T^n]]\) contains \((T^{n-1},T^n)\),
    Nakayama's lemma \cite[Tag 00DV (6)]{stacks-project} implies that there is a surjection of locally free modules of rank one, \(\phi_n \colon R[[T^{n-1},T^n]]\to L|_{V_R}\), which consequently is an isomorphism, finishing the proof of the claim that \(K(R)=\Pic_{X_n}(R)\).

    By the Beauville-Lazlo theorem applied to the Cartier divisor \((T^{n-1})\) (cf \cite{Beauvill-Laszlo} and \cite[Theorem 2.12.1]{Beilinson2005QUANTIZATIONOH}), 
    \[K(R)\cong R[T^{-1}]^{\times}\setminus R((T))^{\times}/R[[T^n,T^{n-1}]]^{\times}\]
    as sets, since \((T^{n})^{n-1}=(T^{n-1})^n\) implies that the \(T^{n-1}\)-adic completion of \(R[T^p,T^{p-1}]\) is isomorphic to the \((T^{n},T^{n-1})\)-adic completion.
    Considering that tensor product of line bundles corresponds to multiplication of any choice of cocycles, this is also an isomorphism of abelian groups. 
    The fact that \(\Pic_{X_n}^{0}(R)\) is the set of classes associated to degree zero Laurent series follows for example by the degree being preserved by pullback along the normalization and that each coset of \(\Pic_{\mathbb{P}^1}(R)=R[T^{-1}]^{\times}\setminus R((T))^{\times}/R[[T]]^{\times}\) has a unique representative of the form \(T^{n}\).

    Assume now that \(n=p=\mathrm{char}\text{ }k\), we want to 
    show that \(\Pic_{X_p}^0\cong \mathbb{G}_a^{\frac{(p-1)(p-2)}{2}}\).
    First, note that as \(X_p\) is a curve, its coherent cohomology (with any coefficients) vanishes in degrees greater than one, hence by \cite[Proposition 8.4.2]{NeronModels}, \(\Pic_{X_p}\) is formally smooth. 
    Moreover, by \cite[Corollary 4.18.3]{Kleiman}, \(\Pic_{X_p}^0\) is a quasi-projective scheme, in particular of finite type, and by the above smooth. In particular, to compute its group scheme structure, it is enough to compute the underlying group variety, which we do in the sequel. We mainly use this to reduce to computing \(\Pic_{X_p}(R)\) for a integral domains \(R\), in particular, \[
    R[T^{-1}]^{\times}=R^{\times} \qquad \text{and} \qquad \Pic_{X_p}(R)\cong R^{\times}\setminus R((T))/R[[T^{p-1},T^{p}]]^\times.
    \]

    To give a concrete description of this double coset, we need to describe \(R[[T^{p},T^{p-1}]]^{\times}\) more explicitly. 
    Thus, we should understand in which degrees there can be non-zero coefficients. To do this, note that an element of \(R[[T^p,T^{p-1}]]^{\times}\) is an invertible power series in \(T\) with potentially non-zero coefficients only in degrees \(d=kp+l(p-1)\) for \(k,l\) non-negative integers, that is in degrees \begin{align}
    \label{equation gaps}
    0,p-1,p,2p-2,2p-1,2p,3p-3,3p-2,3p-1,3p,\ldots,(p-2)p\end{align}
    and in all degrees higher, or equivalently elements appearing in the numerical semi-group generated by \(p,p-1\). Indeed, one sees this by noting that the smallest solution to \(yp+1=x(p-1)\) is \(x=(p-1)\) and \(y=p-2\). In particular, \(\Pic_{X_p}(R)\) is isomorphic to \(\coprod_{n\in \mathbb{Z}} R^{N}\) where \(N\) is the cardinality of \(\mathbb{N}\setminus \langle p,p-1\rangle\) and the index \(n\) ranges over the degree of elements in \(k((T))\) equivalently the degree of the line bundle associated to the cosets.

    Hence, \(\Pic_{X_p}^0\cong \mathbb{A}^N\) as a scheme, and by a theorem of Lazard, \cite[§4, 4.1(Thm. of Lazard)]{DemazureGabriel}, it is unipotent as a group scheme. 
    Note that \(\Pic_{X_p}\) is \(p\)-torsion. Indeed, as the group operation on \(\Pic_{X_p}\) can be described via multiplication of cocycles, and \[f(T)^p=f^{(p)}(T^{p})\in k[[T^{p-1},T^{p}]]\] for any \(f\) of degree \(0\) where \(f^{(p)}(T)=\sum a_i^pT^i\) for \(f(T)=\sum a_iT^i\).
    Since it is also commutative, it follows that it is isomorphic to \(\mathbb{G}_a^N\) as a group scheme; see \cite[Lemma B.1.10]{CGP}. Counting the gaps in the sequence (\ref{equation gaps}) gives the claimed identity \(N=(p-1)(p-2)/2\). 

    \smallskip

    Now, consider \(X_3\) over an algebraically closed base field of arbitrary characteristic. By the above dimension counts (which are independent of characteristic), \(\Pic_{X_3}\cong{\mathbb{A}^1}\) as a scheme. Since there is only one group scheme structure on \(\mathbb{A}^1\) over an algebraically closed field, namely \(\mathbb{G}_a\), we are done.

    Finally, to see that the decomposition is \(\aalpha_p\)-invariant in the characteristic \(p\) case, recall that the coaction in both charts is \(T^{-1}\mapsto T^{-1}+sT^{-p}\), \(T^{p}\mapsto T^p\) and \(T^{p-1}\mapsto T^{p-1}+s\) under the identification \(\aalpha_p=\Spec k[s]/(s^p)\); this action makes the inclusions \begin{align} 
    \label{inclusions R} 
    R[T^{-1}], \, R[[T^p,T^{p-1}]]\subset R((T))\end{align}
    \(\aalpha_p\)-equivariant. Note that this defines a continuous coaction on \(R[[T^p,T^{p-1}]]\) as well as on \(R((T))\), since \(s\) is nilpotent and the action fixes \(T^p\). Indeed, this means that the action on \(R((T))\) is via \(R((T^p))\)-module homomorphisms, and \(R((T))\) is a finite free \(R((T^p))\) module. Analogously, one gets continuity on \(R[[T^p,T^{p-1}]]\), using its \(R[[T^p]]\)-module structure. The fact that this corresponds to the natural pullback action on \(\Pic_{X_p}\) follows from considering what pulling back cocycles along the action map corresponds to on coordinate rings. Namely, let \(\mathcal{L}\in \Pic_{X_p}(R)\) and choose \(f\in R((T))^{\times}\) as a representative cocycle for the corresponding coset. This corresponds to a morphism \(\Spec R((T))\to \mathbb{G}_m\) and precomposing with the action \(a:\aalpha_p\times \Spec R((T))\to \Spec R((T))\) corresponds to the coaction applied to \(f\). 
    Since the inclusions in (\ref{inclusions R}) 
    are $\aalpha_p$-equivariant, the cocycle \[
    f\circ a:\aalpha_p\times \Spec R((T))\to \mathbb{G}_m
    \]
    on \(\aalpha_p\times X_p\) represents the pulled back line bundle \(a^*\mathcal{L}\). The corresponding class 
    \[
    [f \circ a \circ \beta] = [f(T(1-\beta T^{1-p}+\beta ^2T^{-2(p-1)}+\ldots +(-\beta)^{p-1}T^{-(p-1)^2}))]
    \] 
    is functorial in \(R\) and \(\beta:\Spec R\to \aalpha_p\), as it is defined by \(R\)-algebra homomorphisms hence by the Yoneda lemma it determines the action \(\aalpha_p\times \Pic_{X_p}\to \Pic_{X_p}\).
    This finishes the proof.
\end{proof}

\subsection{Quasi-bielliptic, characteristic two and three} Here $X = C$ is the rational curve having a ordinary cusp with a double point, namely the compactification of $(x^3-y^2=0)$ with a point at infinity. By \Cref{lem: picardXp}, the Picard group of $X$ is isomorphic to the additive group. More precisely, for any $k$-algebra $R$, there is a canonical identification $\Pic_X^0(R) \simeq \mathbf{G}_a(R)$, which assigns to the variable $c \in R$ the transition function $1+cs \in R((s))^\times$.  We denote as $L_c$ the class of the line bundle associated to the canonical representative $1+cs$. Notice that this holds in any characteristic.

We shall now exploit the explicit description of $\Pic_X^0$, in order to  compute the infinitesimal quotient of the Picard scheme for a quasi-bielliptic surface $S$. We follow here the classical reference \cite[Propositions 6 and 7]{BMIII}, where such surfaces are completely described in terms of contracted product. For a more recent exposition (from the point of view of group schemes) we refer to \cite[Subsections 5.2.1 and 5.2.2]{Martin}. 

\begin{proposition}
\label{prop: picinf quasi bielliptic}
    Let $S$ be a quasi-bielliptic surface. Then either its Picard scheme is smooth, or its infinitesimal quotient $(\underline{\mathrm{Pic}}^0_S)_{\mathrm{inf}}$ is isomorphic to one among the following:
     \begin{enumerate}
              \item $\boldsymbol{\alpha}_3$ for $p=3$;
          \item either $\boldsymbol{\alpha}_2$ or $\boldsymbol{\alpha}_4$ for $p=2$.
     \end{enumerate}
\end{proposition}

In a late stage of this work, the authors were made aware of the PhD thesis \cite{Kroon}, where quasi-bielliptic surfaces are studied over arbitrary ground fields of positive characteristic. In particular, Section 7.1 computes the quotient $\Pic^\tau_S/(\Pic^0_S)_{\mathrm{red}}$, where $(-)^\tau$ denotes the torsion part, whose connected component of the identity is \((\Pic_S^{0})_{\mathrm{inf}}\). A common tool with the present work is Jensen's results on Picard schemes of quotients (\Cref{thm:Jensen}); besides this, while Kroon's method goes through the study of Bagnera–de Franchis covers of \(S\), we instead rely on the description of \(S\) via \(G\)-normal curves. 
Concerning the statements, \Cref{prop: picinf quasi bielliptic} above has some overlap with \cite[Table 7.1]{Kroon}. In particular, we fill the conjectural entries $(a2)$ in the supersingular quasi-bielliptic case and $(c2)$ in the classical quasi-bielliptic case of the table.

\begin{remark}\label{lem: trivial action on dt}
    In all the actions \(G\times X\to X\) on the cuspidal rational curve we will consider below, we have \(X/G\cong\mathbb{P}^1\). By \Cref{rational quotient case}, letting \(G = G^{lin} \times G^{unip}\) be the factorization of \(G\) into a diagonalizable and a unipotent factor, this means that \(\Pic_{S}^0\) is smooth if and only if \(G^{lin}\) acts non-trivially on \(H^0(X,\omega_X)\). To this end we give a concrete description of a basis of this \(1\)-dimensional \(G\)-representation. 
    As in the previous discussions, let $t$ be the affine parameter for which the cusp is set at infinity, and $\nu \colon \mathbb P^1_t \rightarrow X$ the normalization, given by $x= 1/t^2$ and $y = 1/t^3$. 
 The dualizing sheaf is uniquely determined by its restriction to the smooth part of $X$, so we can restrict to the regular locus of $F(x,y) = y^2-x^3=0$ and compute
\[
\sigma = \frac{\mathrm d x}{F_y} = \frac{\mathrm d x}{2y}, \ (p \neq 2) \qquad \text{or} \qquad \sigma = \frac{-\mathrm d y}{F_x} = \frac{\mathrm d y}{3x^2}, \ (p \neq 3). 
\]
In particular, pulling back and substituting gives $\nu^* \sigma = -\mathrm d t$ as a basis of $H^0(X,\omega_X)$.
\end{remark}

Let us go through the Bombieri-Mumford list and treat the cases of characteristic three and two separately.

\begin{lemma}
    Let $p=3$. Then the only quasi-bielliptic generalized Igusa surface corresponds to $G = \aalpha_3$ of maps $t \mapsto t+bt^3$, with $b^3=0.$ This corresponds to $(d)$ in \cite[page 214]{BMIII}.
\end{lemma}

 \begin{proof}
 The cases $(a)$ and $(b)$ from the Bombieri-Mumford list satisfy $G^{unip} = 1$, so we can exclude them by \Cref{thm:diag_reduced_pic}.
 Next, in case $(c)$ we have $G \supsetneq G^{lin} = \boldsymbol{\mu}_3$, acting via  \[
     t \mapsto at + (1-a)t^3, \qquad a^3= 1,
     \]
     and in case $(e)$, $G \supsetneq G^{lin} = \boldsymbol{\mu}_2$ acting via $t \mapsto \pm t$. In both situations, the $G^{lin}$-action (by pullback) is via a nontrivial character on the differential $\mathrm d t$, so we exclude them by \Cref{lem: trivial action on dt}.  Finally, case $(f)$ does not occur, since $G \simeq \aalpha_9$ is not a subgroup scheme of an elliptic curve, as noted in \cite[Remark 5.12]{Martin}.
 \end{proof}
 
 Now, the gluing datum $t$ on $X$ is stable by the $G$-action, and the Beauville-Lazlo decomposition is compatible with taking pull-backs. Thus, the corresponding action of $G=\boldsymbol{\alpha}_3$ on the coordinate $t$ is given by
\begin{align}
\label{quasihyperelliptic wild char3}
b \cdot t = t+bt^3,
\end{align}
where $b^3=0$ and where the coordinate $t$ is as above. 
 Let $s=1/t$ where defined, then we consider affine coordinates
    \[
    (x,y) = \left(\frac{1}{t^2} = s^2, \frac{1}{t^3}=s^3\right)
    \]
    Let us check that the action (\ref{quasihyperelliptic wild char3}) corresponds to
    \[
    a \cdot (x,y) = (x+a,y), \quad a \in \boldsymbol{\alpha}_3.
    \]
    Indeed, 
    \[
    a \cdot y = a \cdot \frac{1}{t^3} = \frac{1}{a \cdot t^3} = \frac{1}{(t+at^3)^3} = \frac{1}{t^3} = y
    \] since $a^3=0$, and moreover
    \[
    a \cdot x = a \cdot \frac{1}{t^2} = \frac{1}{(t+at^3)^2} = \frac{1}{t^2}+a.
    \]
    To see this last equality, it suffices to factor out $t^2$ from then denominator, and then doing a partial fraction decomposition, using that $p=3=0$ and that $a^3= 0$.  
    Next, since $G= \boldsymbol{\alpha}_3$ is infinitesimal and moreover the $G$-action on $X$ has a fixed point, the assumptions of \Cref{cor:Jensen} are satisfied. So in order to understand what the infinitesimal part of $\Pic_S^0$ is, it suffices to compute the invariant Picard scheme $(\Pic_X^0)^{\boldsymbol{\alpha}_3}$.
    
    \begin{lemma}
     Let $p=3$ and $G = \aalpha_3$ of maps $t \mapsto t+bt^3$, with $b^3=0,$ acting on the cuspidal cubic $X$. The induced $G$-action on $\Pic_X^0 \simeq \mathbb G_a$ is given by 
     \[
     a \cdot L_c = L_{c-ac^3}.
     \]
     In particular, $(\Pic_S^0)_{\mathrm{inf}}$ is isomorphic to $\boldsymbol{\alpha}_3$. This yields (i) in \Cref{prop: picinf quasi bielliptic}.
    \end{lemma}

   \begin{proof}
       It actually suffices to work over the universal test ring $A = k[a,c]/(a^3)$, the ring of functions of $G \times \mathbb G_a$. Then by functoriality the same polynomial law will give the action on any ring $R$ as in the statement. Thus we look for an explicit form of the map
       \[
       a \cdot (-) \colon \Pic_X^0(A) \longrightarrow \Pic_X^0(A), \qquad c \longmapsto a\cdot c = c^\prime
       \]
       for the universal $a \in A$ with $a^3=0$. Consider the line bundle $L_c$ associated to some $c \in \mathbb{G}_a$. Then the action (on the local chart for which the variable $s$ is invertible) is
       \[
       s^\prime := a \cdot s = a \cdot \frac{1}{t} = \frac{1}{t+at^3} = \frac{s^3}{s^2+a}
       \]
       Recall that we are working in the double quotient
       \[
       A[s^{-1}]^\times\backslash A((s))^\times /  A[[s^2,s^3]]^\times,
       \]hence the induced action on the line bundle $L_c$ is given by
       \[
       a \cdot (1+cs) = 1+cs^\prime = 1+c\left(\frac{s^3}{s^2+a}\right) = (a+s^2+cs^3)(s^2+a)^{-1} \sim s^{-2}(a+s^2+cs^3) = 1+cs + as^{-2}
       \]
       because $(s^2+a)^{-1}s^2 = 1-as^{-2}+a^2s^{-4} \in A[s^{-1}]^\times$. It remains to compute some $c^\prime \in A$ such that \[
       1+cs +as^{-2} \sim 1+c^\prime s.
       \]
       Let us do this explicitly: let 
       \[
       w(s) := 1+(ac+a^2c^3)s^{-1} - (a+a^2c^2)s^{-2} +a^2cs^{-3} +a^2s^{-4}.
       \]
       Since $a$ is nilpotent, we have $w \in A[s^{-1}]^\times$,
       so that 
       \[
       1+cs +as^{-2} \sim w(s)(1+cs +as^{-2}) = 1+cs+ac^2+a^2c^4,
       \]
       where the last equality is a direct computation using $a^3=0$ and $3=0$. Then we can multiply the result by $(1-ac^2) \in A^\times$ and this gives
       \begin{align*}
       1+cs +as^{-2} & \sim w(s)(1+cs +as^{-2}) \sim (1-ac^2)w(s)(1+cs +as^{-2}) \\ & = (1-ac^2) ( 1+cs+ac^2+a^2c^4) = 1+(c-ac^3)s
       \end{align*}
       as claimed. Finally, the isomorphism $(\Pic_X^0)_{\mathrm{inf}} \simeq \boldsymbol{\alpha}_3$ follows from computing the kernel of $c \mapsto (c-c^3a)$, which is clearly defined by the ideal $(c^3)$.
       \end{proof}

   \begin{remark}
       It is a natural question to extend the above to the Hilario-Schröer examples. Here are some first steps toward such a computation. Let $S_p = E \times^{\boldsymbol{\alpha}_p} X_p$ be the surface defined in \Cref{subsection:generalized_Igusa_alpha_p}. As proven in \Cref{prop:E}, these have singular Picard scheme when $p\geq 3$, and the case $p=3$ recovers the quasi-bielliptic. Moreover, as in the cusp case above, the $G$-action satisfies the assumptions of \Cref{cor:Jensen}, thus \[
        (\Pic_{S_p}^0)_{\mathrm{inf}} \simeq ((\Pic_{X_p}^0)^G)_{\mathrm{inf}}.
       \]
       So now we want to understand what is the subgroup $(\Pic_{X_p}^0)^G$ and we can focus on the $G$-action on the curve $X_p$. To capture all the information over arbitrary $k$-algebras, it suffices here to work over the universal ring $A=k[a]/(a^p)$. Following (\ref{eq:coaction_Xp}), the co-action of $\boldsymbol{\alpha}_p$ on the variable $T$ is given by 
       \[
       T^\prime := a \cdot T = \frac{1}{a\cdot T^{-1}} = \frac{T^p}{T^{p-1}+a} = T(1+aT^{1-p})^{-1} = T \sum_{m=0}^{p-1}(-aT^{1-p})^m.
       \]
       By \Cref{lem: picardXp}, the connected component of the Picard scheme of the curve $X_p$ is isomorphic to $\mathbb G_a^N$, where $N=(p-1)(p-2)/2$.
       Now we have a natural description: each $c = (c_1,\ldots,c_N) \in \mathbb G_a^N \simeq \Pic^0_{X_p}$ corresponds to the line bundle $L_c$ on $X_c$, with transition function 
\[
\sigma_c(T) = 1+ \sum_{i=1}^N c_iT^i.
\]
To determine an explicit form for the $\boldsymbol{\alpha}_p$-action then means to find $c^\prime \in \mathbb G_a^N$ such that 
\[
a \cdot \sigma_c(T) = \sigma_c(T^\prime) = 1+\sum_{i=1}^N  c_i(T^\prime)^i = 1 + \sum_{i=1}^N  c_i (1+aT^{1-p})^{-i} T^i \sim 1 + \sum_{j=1}^N  c_j^\prime T^j.
\]
in the double quotient. Notice that here, the analogous reduction to integral test-rings as the one in the proof of \Cref{lem: picardXp} cannot be applied. Indeed, since the acting group is the non-reduced \(G= \boldsymbol{\alpha}_p\), the universal ring contains a nilpotent element \(a\). Thus, negative powers of \(T\) naturally occur in the computation. 
\end{remark}

\begin{lemma}
\label{lem: cases f and h}
Let $p=2$. The quasi-bielliptic surfaces having singular Picard scheme are exactly the contracted products $E \times^G X$, where $X$ is the cuspidal cubic and
\begin{itemize}
    \item $G = \aalpha_2$ of maps $t \mapsto t+\lambda at^2+at^4,$ with $a^2=0$;
    \item $G =I_{1,1}$ of maps $t \mapsto t+a+ a^2t^2+a^2t^4$, with $a^4=0$.
\end{itemize}
These correspond respectively to $(f)$ and $(h)$ in \cite[page 214]{BMIII}.
\end{lemma}

\begin{proof}
Looking at the Bombieri-Mumford table, cases $(a)$, $(b)$ and $(d)$ satisfy $G^{unip} = 1$, so we can exclude them by \Cref{thm:diag_reduced_pic}. 
Next, in case $(c)$ we have $G \supsetneq G^{lin} = \boldsymbol{\mu}_2$, acting via \[
t \mapsto at+\lambda(a+1)t^2+(a+1)t^4, \qquad a^2=1;
\]
then in case $(e)$ we have $G \supsetneq G^{lin} = \boldsymbol{\mu}_4$, acting via
\[
t \mapsto bt+(b+b^2)t^2 +(1+b^2)t^4, \qquad b^4=1;
\]
and finally in case $(g)$ we have $G \supsetneq G^{lin} = \boldsymbol{\mu}_3$, acting via $t \mapsto ct^3$ with $c^3=1$. In all of the three last situations, the $G^{lin}$-action (by pullback) is via a nontrivial character on the differential $\mathrm d t$, and we conclude again by \Cref{lem: trivial action on dt}.
\end{proof}

Let us place ourselves under the following assumption: $G \simeq \boldsymbol{\alpha}_2$, acting by
\begin{align}
\label{quasihyperelliptic wild char2}
a \cdot t = t + \lambda at^2 + at^4
\end{align}
for some $\lambda \in k$, where $a^2 = 0$ is the formal variable for $\boldsymbol{\alpha}_2$. The action translates to
\[
a\cdot (x,y) = (x, y+\lambda ax + a), \qquad a \in \boldsymbol{\alpha}_2.
\]
Indeed, we have 
\[
a \cdot x = a\cdot \frac{1}{t^2} = \frac{1}{a\cdot t^2} = \frac{1}{(t+\lambda at^2+at^4)^2} = \frac{1}{t^2} = x,
\]
since $a^2 = 0$, and analogously
\[
a \cdot y = a \cdot \frac{1}{t^3} = a \cdot \frac{1}{t^2} \, a \cdot \frac{1}{t} = \frac{1}{t^2} \frac{1}{a\cdot t} = \frac{1}{t^2(t+\lambda a t^2 + at^4)} = \frac{t+\lambda a t^2 + at^4}{t^4} = y + \lambda a x+ a,
\]
as claimed. Since $G$ is infinitesimal and moreover the $G$-action on $X$ has a fixed point, the assumptions of \Cref{cor:Jensen} are again satisfied. Thus, it suffices to understand the invariant Picard scheme $(\Pic_X^0)^{\boldsymbol{\alpha}_2}$.

\begin{lemma}
Let $p=2$ and $G = \aalpha_2$ acting on the cuspidal cubic $X$ as in (\ref{quasihyperelliptic wild char2}). The induced $G$-action on $\Pic_X^0 \simeq \mathbb G_a$ is given by 
     \begin{align}
     \label{lem: case f cusp}
     a \cdot L_c = L_{c+\lambda ac^2+ac^4}.
     \end{align}
     In particular, $(\Pic_S^0)_{\mathrm{inf}}$ is isomorphic to $\boldsymbol{\alpha}_4$ when $\lambda = 0$, and to $\boldsymbol{\alpha}_2$ otherwise. This yields (ii) in \Cref{prop: picinf quasi bielliptic}. 
    
    \end{lemma}

\begin{proof}
As in the characteristic three case, it actually suffices to work over the universal test ring $A= k[a]/(a^2)$. Consider the line bundle $L_c$ associated to some $c \in \mathbb G_a(A) = A$. The action on the variable $s$ is
\[
s^\prime := a \cdot s = a \cdot \frac{1}{t} = \frac{1}{t+\lambda at^2+at^4}= s (1+\lambda as^{-1} + as^{-3})
\]
Now, what remains to be computed is $c^\prime$ such that 
\[
1+cs^\prime = 1+\lambda ac + cs + acs^{-2} \sim 1+c^\prime s
\]
Multiplying the left hand side by the unit \[
w(s) = 1+a (\lambda c + c^3 +c^2s^{-1}+cs^{-2}) \in A[s^{-1}]^\times\] gives
\[
w(s)(1+cs^\prime) = 1+s (c+\lambda ac^2+ac^4), 
\]
so $c^\prime = c+\lambda ac^2+ac^4$ as claimed. Finally, the subgroup scheme of fixed point is given by the kernel of the additive polynomial $c^4+\lambda c^2$. If $\lambda =0$ we get the second Frobenius kernel of the additive group. If $\lambda \neq 0$ the equation factors as $c^2(c+\lambda)^2=0$ and its infinitesimal part is thus isomorphic to $\boldsymbol{\alpha}_2$.
\end{proof}

The very last situation that is left to consider is the second item from \Cref{lem: cases f and h}, where the acting group is the $2$-torsion of a supersingular elliptic curve, which we denote as $I_{1,1}$ (in \cite{Martin} this subgroup scheme is called $M_2$).

\begin{lemma}
Let $p=2$ and $G = I_{1,1}$ of maps $t \mapsto t + a+ a^2t^2+a^2t^4$, with $a^4=0,$ acting on the cuspidal cubic $X$. Then $(\Pic_S^0)_{\mathrm{inf}}$ is isomorphic to $\boldsymbol{\alpha}_2$. This is covered by case (ii) in \Cref{prop: picinf quasi bielliptic}. 
\end{lemma}

\begin{proof}
Notice that the $G$-action on $X$ has no fixed point here, so that \Cref{cor:Jensen} cannot be applied directly. Let us proceed by a dévissage argument: let $H \simeq \aalpha_2 \subset G$ of maps $t\mapsto t+a$ ($a^2=0$), with quotient $K := G/H \simeq \aalpha_2$ acting on the intermediate contracted product $W:= E \times^H X$. We are interested in the surface $S$ in the following
\[
E \times X \stackrel{/H}{\longrightarrow} W \stackrel{/K}{\longrightarrow}  S := E \times^G X = W/K.
\]
Let us look at the $H$-action on the affine chart of $X$, for which the cusp is at the origin; there we have
\[
x = \frac{1}{t^2} \longmapsto a \cdot x = \frac{1}{(t+a)^2} =x, \qquad y = \frac{1}{t^3} \longmapsto a \cdot y = \frac{1}{(t+a)^3} = \frac{1}{t^2(t+a)} = \frac{t+a}{t^4} = y + ax^2
\]
Since the only invariant functions are those independent of the variable $y$, as a quotient we get $X/\aalpha_2 = \mathbb P^1$. The cusp, which we denote as $\{t_\infty\}$ is a fixed point by the action. Thus the map $E \longrightarrow E \times \{t_\infty\}$ is $H$-equivariant and yields a section $\sigma \colon E/H \to W$. This way, the sequence
\[
0 \longrightarrow \Pic^0_X \longrightarrow \Pic^0_W \longrightarrow \Pic^0_{E/H}\longrightarrow 0
\]
splits. Now, the $K$-action on the direct product $\Pic_W^0 = \Pic^0_{E/H} \times \Pic^0_X$ is by definition an homomorphism \[
K\to \mathrm{Aut}_{\mathrm{gp}}(\Pic^0_{E/H} \times \Pic^0_X)\cong \mathrm{Aut}_{\mathrm{gp}}(\Pic^0_{E/H})\times \mathrm{Aut}_{\mathrm{gp}}(\Pic^0_X),\] 
where the last equality comes from the fact that there are no nontrivial group homomorphisms between an elliptic curve and the additive group. In other words, the $K$-action is diagonal. On the first factor, it is the trivial one, since it is via translations on $E/H$ coming from translations on $E$. On the second factor, we claim that the $K$-action identifies with
\[
K \simeq \aalpha_2 \curvearrowright \Pic^0_X \simeq \mathbb G_a = \Spec( k[c]), \qquad c \longmapsto b \cdot c = c+bc^2+bc^4 \quad (b^2=0).
\]
Indeed, let us check it over the universal ring $A= k[a]/(a^4)$, using the double coset description of \Cref{lem: picardXp}. The action on $s=1/t$ is given by
\[
s \longmapsto s^\prime = a \cdot s = \frac{1}{a\cdot t}= \frac{1}{t+a+a^2t^2+a^2t^4}= s+as^2+a^2+a^2s^{-2}+a^2s^3+a^3s^4.
\]
Notice that we can work modulo $s^2A[[s]] \subset A[[s^2,s^3]]^\times$, and that the $H$-action on $\Pic_X^0$ is trivial. In particular, on $\Pic_X^0(A) = \mathbb G_a(A) = A[c]$, we have that $c \mapsto c^\prime = a\cdot c$ where
\[
1+c^\prime s = 1+c(a\cdot s) \sim 1+cs +ca^2+ca^2s^{-2},
\]
which, multiplying by $w = 1+a^2c^2s^{-1} +a^2cs^{-2} \in A[s^{-1}]^\times$ and then taking modulo $s^2A[[s]]$ again, gives
\[
1+c^\prime s \sim 1+a^2(c+c^3)+cs = (1+a^2(c+c^3))(1+(c+a^2c^2+a^2c^4)s).
\]
Thus, the $K$-action on the variable $c$ is obtained by setting $b=a^2$ which gives $c^\prime = c+bc^2+bc^4$ as wanted.
Next, we claim that
\begin{align}
\label{last equality}
(\Pic_S^0)_{\mathrm{inf}} = ((\Pic_W^0)^K)_{\mathrm{inf}}.
\end{align}
This suffices to conclude our proof; indeed,
\[
((\Pic_W^\tau)^K)_{\mathrm{inf}} =  ((\Pic^0_X)^K)_{\mathrm{inf}}  = ({\mathbb G}_{a}^{\aalpha_2})_{\mathrm{inf}}.
\]
Notice that the $K$-action on $\mathbb G_a$ is exactly as in (\ref{lem: case f cusp}) with $\lambda = 1$, and this implies \[
(\Pic_S^0)_{\mathrm{inf}}  = ({\mathbb G}_{a}^{\aalpha_2})_{\mathrm{inf}} \simeq \aalpha_2.\] By applying \Cref{coro:B} to the $G$-action on $X$, we know that the Picard scheme of $S$ must be singular. So the inclusion is an equality, and we conclude that $(\Pic_S^0)_{\mathrm{inf}} = \aalpha_2$. As a last step, let us prove (\ref{last equality}); from \Cref{thm:Jensen}, we get the short exact sequence
\[
0 \longrightarrow D(K) \longrightarrow \Pic_S^\tau \longrightarrow (\Pic_W^\tau)^K \longrightarrow 0
\] 
where the Cartier dual is $D(K) \simeq \aalpha_2$. 
    Let $p_W\colon W \to E/H$ and $h\colon S\to E/G$ be the natural projections. Analogously to the proof of \Cref{cor:Jensen}, we have the following commutative diagram
    \[
    \begin{tikzcd}
        0 \ar[r] & D(K) \ar[r] & \Pic_S^\tau \ar[r] &  (\Pic_{W}^\tau)^K \ar[r] & 0 \\
        0 \ar[r] & D(K) \ar[r]\ar[u,equal] & \Pic_{E/G}^\tau \ar[r] \ar[u,"h^*"] &  (\Pic_{E/H}^\tau)^K \ar[r] \ar[u,"p_W^*"] & 0 
    \end{tikzcd}
    \]
   where the square on the right hand side is cartesian and the two horizontal lines are short exact sequences.

    By the commutativity of the diagram above, it follows that 
    \[D(K)=\ker\left(\Pic_S^{\tau}\to (\Pic^{\tau}_{{W}})^K\right)\subset h^*(\Pic_{E/G}^{\tau})\subset (\Pic_{S}^{\tau})_{\mathrm{red}}\]
    where the last inclusion comes from \(\Pic_{E/G}^{\tau}\) being an abelian variety. Now, the restriction of \(\Pic_{S}^{\tau}\to (\Pic_{W}^{\tau})^K\) to \((\Pic_{S}^{\tau})_{\mathrm{red}}\) factors through \(((\Pic_{W}^\tau)^K)_{\mathrm{red}}\), by the universal property of underlying reduced subschemes. In particular, \(\Pic_S^\tau\to (\Pic^\tau_{{W}})^K\) induces a surjective morphism \((\Pic^\tau_S)_{\mathrm{inf}}\to ((\Pic_{W}^\tau)^K)_{\mathrm{inf}}\), whose kernel is trivial, and hence it is an isomorphism.
\end{proof}

\subsection{Bielliptic, characteristic two and three} 
We consider here generalized Igusa surfaces among bielliptic ones, which only exist in small characteristics. Such surfaces are described completely as contracted products of the form $E \times^G X$ in \cite[page 37]{BMII}, where $E$ is an ordinary elliptic curve and $X = E_0$ is a second elliptic curve.

 \begin{proposition}
 \label{prop: picinf bielliptic}
     Let $S$ be a bielliptic surface. Then either its Picard scheme is smooth, or its infinitesimal quotient $(\underline{\mathrm{Pic}}^0_S)_{\mathrm{inf}}$ is isomorphic to one among the following:
     \begin{enumerate}
          \item $\boldsymbol{\mu}_2$ for $p=2$;
         \item the $2$-torsion of a supersingular elliptic curve, for $p = 2$;
         \item $\boldsymbol{\alpha}_p$, for $p \in \{2,3\}$.
     \end{enumerate}
 \end{proposition}

 \begin{proof}
 Reinterpreting what is said in the proof of \cite[Theorem 1.2] {Liedtke_note} in an explicit way gives a complete list of bielliptic surfaces which are generalized Igusa; notice that in all of the surfaces listed below, either the group $G$ is unipotent, or there is a diagonalizable factor acting by translation; thus the dualizing sheaf of $E_0$ has non-trivial invariant sections (see \Cref{lem:invariants}). Moreover, either $G$ satisfies one of the assumptions in \Cref{thm:Jensen}, or it is a direct product of two factors satisfying such assumptions. If moreover the $G$-action on $X$ has a fixed point, by \Cref{cor:Jensen} we can compute the infinitesimal quotient as
   \[
    (\underline{\mathrm{Pic}}^0_S)_{\mathrm{inf}} = (\underline{\mathrm{Pic}}_{E_0}^0)^G_{\mathrm{inf}} = ((E_0^\vee)^G)_{\mathrm{inf}},
    \]
where $(-)^\vee$ denotes the dual elliptic curve.

When $p=3$, the only generalized Igusa surface obtained this way is given by $G=\mathbb{Z}/3\mathbb{Z}$ acting on $E_0$ with an automorphism of order three that is not a translation; the condition is $j(E_0) =0$ which means that the curve is supersingular. The $G$-action on the dual curve is given by the dual automorphism; so we get $(\underline{\mathrm{Pic}}^0_S)_{\mathrm{inf}} = (E_0^\vee)^G$. To compute the latter, since $E_0$ is unique up to isomorphism, we may choose an embedding into the projective plane $\mathbb P^2 = \mathbb P^2_{x,y,z}$ as $E_0^\vee = (zy^2=x^3+xz^2)$. As described in \cite[Appendix A, Proposition 1.2]{SilvermanEllCurves}, the automorphism group of $E_0$ (up to translations) is finite of order $12$, and the order three element is given by $(x\colon y \colon z) \mapsto (x+zr \colon y \colon z)$, where $r^3+r=0$. Geometrically, the only fixed point is then $(0\colon 1 \colon 0).$ Working locally in the affine coordinates $(x,z)$ around the origin, one sees that the fixed scheme is given by $(z=x^3=0)$. Again using that $E_0$ is supersingular (and thus its only infinitesimal subgroup scheme of length three is the additive Frobenius kernel) we conclude  
\[
((E_0^\vee)^G)_{\mathrm{inf}} = ((E_0^\vee)^G) \simeq \boldsymbol{\alpha}_3,
\]
covering (3) of the statement for $p=3$.

    When $p=2$, the bielliptic surfaces with singular Picard are those in cases $(a)$ and $(c)$ of the list in \cite{BMII}. The first case is the Igusa surface, for which
   \Cref{cor:Jensen} yields $(\Pic_S^0)_{\mathrm{inf}} = (E_0^\vee[2])_{\mathrm{inf}}$. The latter is either equal to $\boldsymbol{\mu}_2$ or to a nontrivial extension of $\boldsymbol{\alpha}_2$ by itself, according to whether $E_0$ is ordinary or not; this fits into (1) and (2) in the statement.
   
Next, one has to consider
    \[
    G := E[2] = (\mathbb{Z}/2\mathbb{Z}) \times \boldsymbol{\mu}_2,
    \]
    where the $\mathbb{Z}/2\mathbb{Z}$-action is the same as in the Igusa case, and the $\boldsymbol{\mu}_2$ acts by translation on both factors. In this case, \Cref{cor:Jensen} does not directly apply, because the $\boldsymbol{\mu}_2$-action on $X= E_0$ (here $E_0$ is an ordinary elliptic curve) has no fixed point. Thus let us denote as $ W = E \times^{\mathbb Z/ 2\mathbb Z} E_0$ the Igusa surface, so that $S = W/\boldsymbol{\mu}_2$. Then
    \Cref{thm:Jensen} yields a short exact sequence
    \[
    0 \rightarrow \mathbb{Z}/2\mathbb{Z}\rightarrow \Pic_S^0\rightarrow (\Pic_W^0)^{\mu_2} \rightarrow 0
    \]
    In particular, taking infinitesimal quotients yields \[
    (\Pic_S^0)_{\mathrm{inf}}  = (\Pic_W^{\boldsymbol{\mu}_2})_{\mathrm{inf}} \subset (\Pic_W^0)_{\mathrm{inf}} = \muu_2,\]
    where the last equality comes from the Igusa case. Since the $G^{lin}= \muu_2$-action on $H^0(E_0,\omega_{E_0})$ is trivial, the Picard scheme must be singular by \Cref{rational quotient case}. 
    This forces $(\Pic_S^0)_{\mathrm{inf}}  = \muu_2$, which fits into (1) 
    in the statement.
    
    Finally, let us consider $G=\mathbb{Z}/4\mathbb{Z}$ acting on the elliptic curve $X$ via an automorphism of order four, with the condition that the $j$-invariant of $E_0$ is equal to $12^3$, so that it is a supersingular elliptic curve. This yields $(\Pic_S^0)_{\mathrm{inf}} = (E_0^\vee)^G$. To compute the latter, we may choose the embedding of the dual curve into the projective plane $\mathbb P^2 = \mathbb P^2_{x,y,z}$ given by $E_0^\vee = (zy^2+z^2y=x^3)$. Analogously to the characteristic three case, we follow \cite[Appendix A, Proposition 1.2]{SilvermanEllCurves}, where the automorphism group of $E_0$ is described up to translations. The action of the nontrivial automorphism of order four is given by $(x \colon y \colon z ) \mapsto (x+z \colon y+x+zt \colon z)$, where $t^2+t=1$. Geometrically, the only fixed point is $(0\colon 1\colon 0)$, so we can consider affine coordinates $(x,z)$ and work around the origin. The action becomes  \[
    (x,z) \mapsto \left( \frac{x+z}{1+x+tz}, \frac{ z}{1+x+zt}\right)
    \]
    and thus the locus of fixed points is
    \[
    z=x^2+tzx, \qquad zx +z^2t = 0.
    \]
    Multiplying the first equation by $z$ and then substituting from the second equation gives $z^2 = x^2z+x^2z =0$, which yields $z^2=0$. Then from the second equation one gets $xz=0$, so either $x=0$ or $z=0$. If $x=0$ then $z=0$ (and we get back only the reduced structure of the point), while the non-reduced structure comes from $z=0$, which implies $x^2=0$. Finally, the fixed point subscheme is given by $(z=x^2=0)$. Since the only infinitesimal subgroup scheme of length two of a supersingular elliptic curve is $\boldsymbol{\alpha}_2$, we conclude
    \[
    ((E_0^\vee)^G)_{\mathrm{inf}} = ((E_0^\vee)^G) \simeq \boldsymbol{\alpha}_2,
    \]
    covering (3) of the statement for $p=2$.
    \end{proof}

\printbibliography

\end{document}